\documentclass[a4paper,10pt]{article}
\usepackage{latexsym,amsthm,psfrag,amsfonts,amsmath,epsfig,comment}

\makeatletter
\@addtoreset{figure}{subsection}
\makeatother


\begin{document}

\thispagestyle{empty}

\newtheorem{lemma}{Lemma}[subsection]
\newtheorem{cor}{Corollary}
\newtheorem{theorem}{Theorem}
\renewcommand{\proofname}{Proof}
\newtheorem{property}{Property}
\newtheorem{prop}{Proposition}[subsection]

\theoremstyle{definition}
\newtheorem{defin}{Definition}
\newtheorem{remark}{Remark}

\def \H{{\mathbb H}}
\def \R{{\mathbb R}}
\def \E{{\mathbb E}}
\def \S{{\mathbb S}}
\def \G{{\mathcal G}}

\def \l{\langle }
\def \r{\rangle }
\def \[{[ }
\def \]{] }
\def \d{D\,}
\def \sign{\text{\,sign\,}}
\def \conv{\text{\,conv\,}}
\def \wt{\widetilde}
\def \wh{\widehat}
\def \a{\alpha}
\def \b{\beta}
\def \dim{\text{\,dim}}

\newcommand{\vn}{\mathop{\mathrm{int}}\nolimits}
\newcommand{\rel}{\mathop{\mathrm{rel\:int}}\nolimits}
\newcommand{\w}{\widetilde }
\renewcommand{\o}{\overline }

\begin{center}
{\huge
Coxeter polytopes with a unique pair of non-intersecting facets.
}

\vspace{20pt}

\begin{tabular}{cc}
{\Large  Anna Felikson}\footnotemark[1] \phantom{asasas} &
\phantom{asasas}
{\Large  Pavel Tumarkin}\footnotemark[2]\\
felikson@mccme.ru \phantom{asasas} &
\phantom{asasas} pasha@mccme.ru \\

\\
\end{tabular}
\footnotetext[1]{Partially supported by grants NSh-5666.2006.1, 
INTAS YSF-06-10000014-5916, and RFBR 07-01-00390-a.}
\footnotetext[2]{Partially supported by grants
MK-6290.2006.1, NSh-5666.2006.1, INTAS YSF-06-10000014-5766, and RFBR 07-01-00390-a.}

\begin{tabular}{c}
Independent University of Moscow, Russia\\
%
%
University of Fribourg, Switzerland
\end{tabular}

\end{center}
\vspace{20pt}

\begin{center}

\parbox{11.5cm}{
\small
{\it Abstract.}
We consider compact hyperbolic Coxeter polytopes whose Coxeter diagram contains a unique dotted edge.
We prove that such a polytope in $d$-dimensional hyperbolic space has at most $d+3$ facets. 
In view of~\cite{Lan},~\cite{K},~\cite{Ess}, and~\cite{d+3}, this implies that   
compact hyperbolic Coxeter  polytopes with a unique pair of non-intersecting facets are completely classified.
They do exist up to dimension $6$ and in dimension $8$ only. 

}

\end{center}

\noindent
\vspace{30pt}


\section{Introduction}

We study  compact Coxeter polytopes in hyperbolic spaces.
Besides the general restriction $d<30$ on the dimension $d$ of the polytope \cite{abs}
and investigation of arithmetic subgroups, 
there are several directions in which some attempts of general classification were
undertaken. One of them is to fix the dimension of polytope. 
Compact hyperbolic Coxeter polytopes of dimensions $2$ and $3$ were completely classified in~\cite{P} and~\cite{Andr}, 
respectively. Another direction is to fix the difference between the number of facets of the polytope and 
its dimension. Simplices were classified in~\cite{Lan}, $d$-dimensional polytopes
with $d+2$ facets were classified in~\cite{K} and~\cite{Ess}, 
$d$-dimensional polytopes with $d+3$ facets were classified 
in~\cite{Ess} and~\cite{d+3}. This paper is devoted to investigation of another direction in classification:
the number of pairs of non-intersecting facets.  

In paper~\cite{nodots} we classified all  compact hyperbolic Coxeter polytopes with mutually intersecting facets.
It turns out that they do exist up to dimension 4 only, and have at most 6 facets.
In this paper we expand the technique developed in~\cite{nodots}  
to investigate compact hyperbolic Coxeter polytopes with exactly one pair of non-intersecting facets.
The paper is devoted to the proof of the following theorem:

\bigskip
\noindent
{\bf Main Theorem.}
{\it
A compact hyperbolic Coxeter $d$-polytope with exactly one pair of non-intersecting facets has at most
$d+3$ facets. In particular, no such polytopes do exist in dimensions $d\ge 9$ and $d=7$. 
}

\bigskip

Clearly, neither simplices nor products of simplices (except prisms) have non-intersecting facets. Therefore, the Main Theorem can be reformulated in the following way.

\noindent
{\bf Corollary.}
{\it
Any compact hyperbolic Coxeter $d$-polytope with exactly one pair of non-intersecting facets is either a prism or a polytope with $d+3$ facets. 
}

\bigskip

The proof is based  on already obtained classifications of polytopes 
of either smaller dimensions or with smaller number of facets,
or with smaller number of pairs of non-intersecting facets.
In fact, the technique we use may lead to the inductive algorithm of classification
of compact hyperbolic polytopes with respect to three directions described above. 
In this context the Main Theorem may be considered as the adjusting of the base of
the tentative algorithm.

The paper is organized as follows: in Section~\ref{tool} we expand the technique
developed in~\cite{nodots}  
to the case of compact hyperbolic Coxeter polytopes with exactly one pair of non-intersecting facets.
In Section~\ref{proof} we prove the Main Theorem moving from smaller dimensions to larger ones (namely, 
up to dimension $12$). Then we finish the proof considering dimensions greater than $12$. In the Appendix 
we reproduce the list of all the compact hyperbolic Coxeter polytopes with exactly one pair of non-intersecting facets.

The paper was written during the authors' stay at the University of Fribourg, Switzerland.
We are grateful to the University for hospitality.


\section{ Technical tools}
\label{tool}

We refer to~\cite[Sections~2 and~3.1]{nodots} for all  essential preliminaries.
Concerning Coxeter polytopes and Coxeter diagrams, we mainly follow~\cite{abs} and~\cite{refl}.
We use the technique of local determinants developed in~\cite{abs}.   
Description of Coxeter facets may be found in~\cite{All}. 
We use standard notation for elliptic and parabolic diagrams (see~\cite{refl}).

\subsection{Notation }
We recall some notation introduced in~\cite{nodots}.

We write {\it $d$-polytope} instead of ``$d$-dimensional polytope'' and {\it $k$-face } instead of ``$k$-dimensional face''.
We say that an edge of Coxeter diagram is {\it multiple} if it is of multiplicity at least two, and an edge  is  
{\it multi-multiple}  if it is of multiplicity at least four. For nodes $x$ and $y$ of a Coxeter diagram $\Sigma$
we write  $[x,y]=m$ if  $x$ is joined with $y$ by an $(m-2)$-tuple edge (or by an edge labeled by $m$).
We write  $[x,y]=\infty$ if $x$ is joined with $y$ by a dotted edge, 
and we write $[x,y]=2$ if the nodes $x$ and $y$ are not joined.

If $\Sigma_1$ and $\Sigma_2$ are subdiagrams of a Coxeter diagram $\Sigma$,
we denote by $\l \Sigma_1, \Sigma_2\r$ a subdiagram of $\Sigma$ spanned by all nodes of
$\Sigma_1$ and $\Sigma_2$.
By  $\Sigma_1\setminus \Sigma_2$ we denote a subdiagram of $\Sigma$ spanned by all nodes of
$\Sigma_1$ that do not belong to $\Sigma_2$.
By $|\Sigma|$ we denote an order of the diagram $\Sigma$.

Given a Coxeter $d$-polytope $P$ we denote by $\Sigma(P)$ the Coxeter diagram of $P$.
If $S_0$ is an elliptic subdiagram of $\Sigma(P)$, we denote by $P(S_0)$
the face defined by this subdiagram (it is  a $(d-|S_0|)$-face
obtained by the intersection of the facets corresponding to the nodes of $S_0$). 
We say that $x\in \Sigma(P)$ is a {\it neighbor} of $S_0$ 
if $x$ attaches to $S_0$ (i.e. $x$ is joined with $S_0$ by at least one edge),
otherwise we say that $x$ is a {\it non-neighbor} of $S_0$.
We say that the neighbor $x$ of $S_0$ is  {\it good} if 
$\l S_0, x\r$ is an elliptic diagram, and
{\it bad} otherwise.
We denote by $\overline S_0$ the sub\-dia\-gram of $\Sigma(P)$
consisting of nodes corresponding to facets of $P(S_0)$. The
diagram $\overline S_0$ is spanned by all good neighbors and all non-neighbors of $S_0$ 
(in other words, $\overline S_0$ is spanned by all vertices of $\Sigma(P)\setminus S_0$ except
bad neighbors of $S_0$). If $P(S_0)$ is a Coxeter polytope, denote
its Coxeter diagram by $\Sigma_{S_0}$.

It is shown in~\cite[Theorem~2.2]{All} that if $S_0$ is an elliptic diagram containing no $A_n$ and $D_5$ components,
then the face $P(S_0)$ is a Coxeter polytope, and its diagram $\Sigma(S_0)$ can be easily
found from the subdiagram $\l S_0,\Sigma_{S_0}\r$.
This fact is the main tool for our induction: if $S_0$ has no good neighbors (this is always the case if $S_0$ 
is of the type $H_4$, $F_4$ or $G_2^{(k)}$, where $k\ge 6$) then $\Sigma_{S_0}=\o S_0$
is a diagram of a Coxeter polytope of lower dimension. 
If the initial polytope has at most one pair of non-intersecting facets, 
then the same is true for $P(S_0)$. So, in assumption that the Main Theorem holds in lower dimensions,   
this implies that $P(S_0)$ is either a simplex, or a triangular prism, or one of 7 Esselmann polytopes,
or one of finitely many $d'$-polytopes with $d'+3$ facets which have diagrams containing at most one dotted edge
(more precisely, in the latter case there are eight 4-polytopes,
a unique polytope 5-polytope  (Fig.~\ref{9_}(c)),
three 6-polytopes (Fig.~\ref{8_multi}),  
a unique 8-polytope (Fig.~\ref{10_}), and no polytopes in dimension 7 and in dimensions greater than 8.

We will also use the following lemmas.

\begin{lemma}
\label{polygon}
Let $P$ be a compact Coxeter $d$-polytope with exactly one pair of non-intersecting facets, and let $S_0\subset \Sigma(P)$ be an elliptic subdiagram.
If $P(S_0)$ is a 2-polytope (i.e. $P(S_0)$ is a polygon) then\\

1) If $\o S_0=\Sigma_{S_0}$ and $\o S_0$ contains no dotted edge, then $\o S_0$ is a Lann\'er diagram of order 3.

2) If $\o S_0$ contains a dotted edge, then $S_0$ has at least one good neighbor.  

\end{lemma}

\begin{proof}
A triangle is the only polygon with mutually intersecting facets, which proves the first statement.
Suppose that $S_0$ has no good neighbors, then $\o S_0=\Sigma_{S_0}$ is a Coxeter diagram of a polygon. Thus,  $\o S_0$ either contains no dotted edges (if $P(S_0)$ is a triangle), or contains at least two. The latter is impossible by the assumptions on $P$, the former contradicts the assumption of the second statement.

\end{proof}

\begin{lemma}
\label{l'}
Suppose that $P$ is a compact Coxeter $d$-polytope with exactly one pair of non-intersecting facets
and at least $d+4$ facets.
Let $\Sigma_1\subset \Sigma(P)$ be a subdiagram of order not greater than $d+2$.
Then

1) There exists a node $x\in \Sigma(P)\setminus \Sigma_1$ such that the subdiagram
$\l x,\Sigma_1\r$ contains no dotted edges. 

2) Suppose in addition that $S\subset \Sigma_1$ is an elliptic diagram  of order $|S|<d$
having less than $d-|S|$ good neighbors and non-neighbors in total in $\Sigma_1$.
Then there exists a node $x\in \Sigma(P)\setminus \Sigma_1$
such that $x$ is not a bad neighbor of $S$ and  the subdiagram
$\l x,\Sigma_1\r$ contains no dotted edges.

3) The statement of the preceding item is also true if $S_1$ has exactly  
$d-|S|$ good neighbors and non-neighbors in total in $\Sigma_1$
and $S_1$ contains an end of the dotted edge.
 
\end{lemma}

\begin{proof}
To prove the first statement, notice that
$\Sigma(P)\setminus \Sigma_1$ contains at least two nodes, at least one of these nodes is not joined with   
$\Sigma_1$ by a dotted edge.
The same consideration works for the second statement: $\o S$ must have at least $d-|S|+1$ nodes, 
so $\Sigma(P)\setminus \Sigma_1$ contains
at least two good neighbors or non-neighbors of $S$.
To prove the third statement, notice that $\Sigma(P)\setminus \Sigma_1$ contains a good neighbor or a non-neighbor of $S$,
which definitely cannot be an end of the dotted edge.

\end{proof}

\subsection{Lists $L(S_0,d)$, $L_1(S_0,d)$ and $L'(\Sigma,C,d)$}

In~\cite[Lemma~3]{nodots} we have proved that if a Coxeter diagram of a polytope contains no dotted edges, then
it contains a subdiagram satisfying some special properties. We have defined a finite list $L(S_0,d)$
of diagrams satisfying these properties.
In this section we slightly change this definition to be applied to the case of diagrams
containing  a unique dotted edge.

We will need the following definitions.

If $\Sigma$ is a Coxeter diagram of a simplicial prism, then the node $x\in \Sigma$ is called a {\it tail }
 if $x$ is an end of the dotted edge and $\Sigma\setminus x$ is a connected diagram.
By a {\it diagram without tail} we mean $\Sigma$ with exactly one of its tails discarded.

We introduce  a partial order ``$\prec$'' on the set of connected elliptic subdiagrams
of maximal order of Lann\'er diagrams and diagrams of simplicial prisms without tail:
\begin{itemize}
\item
 $A_2 \prec B_2\prec G_2^{(k)}$, $k>2$, and $ G_2^{(k)}\prec G_2^{(l)}$ if $k<l$; 
\item
$A_3 \prec B_3 \prec H_3$;  
\item
 $A_4\prec B_4 \prec F_4 \prec H_4$.
\end{itemize}

\medskip
\noindent
{\bf Remark.}
We do not need to introduce a partial order on the diagrams of order 5,
since any diagram of a 5-prism without tail contains connected elliptic diagrams of order 5 of one type only.

\medskip

Suppose that $\Sigma$ is a Lann\'er diagram or a diagram of a simplicial prism without tail.
A connected  elliptic  subdiagram $S\subset \Sigma$ of maximal order is called
{\it maximal in $\Sigma$} if $\Sigma$ contains no connected elliptic subdiagram $S'$ such that $S\prec S'$.
 A connected  elliptic  subdiagram $S\subset \Sigma$ of maximal order  is called {\it next to maximal in $\Sigma$}
if  $\Sigma$ contains a maximal connected elliptic subdiagram $S'$ such that $S\prec S'$ while
 $\Sigma$ contains no connected elliptic subdiagram $S^{''}$ such that $S\prec S^{''} \prec S'$.

\begin{lemma}
\label{l1}
Let $P$ be a compact Coxeter $d$-polytope with a unique pair of non-intersecting
facets, and assume that $P$ has at least $d+4$ facets. Let $S_0$ be a
connected elliptic sub\-dia\-gram of
$\Sigma(P)$ such that \\
{\sc (i)} $|S_0|< d$ and $S_0\ne A_n$, $D_5$.\\
{\sc (ii)} $S_0$ has no good neighbors in $\Sigma(P)$.\\
{\sc (iii)} If $|S_0|\ne 2$, then
$\Sigma(P)$ contains no multi-multiple edges.\\
%
\phantom{{\sc (iii)}} If $|S_0|=2$,
then the edge of $S_0$ has the maximum multiplicity amongst all
edges in $\Sigma(P)$.

Suppose that the Main Theorem holds for any $d_1$-polytope
satisfying $d_1<d$.
Then there exist a sub\-dia\-gram $S_1\subset\Sigma(P)$
and two vertices $y_0,y_1\in\Sigma(P)$ such that the sub\-dia\-gram
$\l S_0,y_1,y_0,S_1\r $ satisfies the following conditions:

\begin{itemize}

\item[$(0)$]
$\l S_0,y_1,y_0,S_1\r $ contains no dotted edges and parabolic subdiagrams;

\item[$(1)$]
    $S_0$ and $S_1$ are elliptic diagrams, $S_0$ is connected, and
    $S_0\ne A_n, D_5$;

\item[$(2)$]
    No vertex of $S_1$ attaches to $S_0$, and 
    $|S_0|+|S_1|=d$;

\item[$(3)$]
    $\l y_0,S_1\r $ is either a Lann\'er diagram, or a diagram of a simplicial prism with a tail discarded,
or one of the diagrams shown in Table~\ref{list} (in the latter case $y_0$ is the marked vertex of the diagram);

\item[$(4)$]
    $\l S_0,y_1\r $ is an indefinite sub\-dia\-gram, and
    $y_1$ is either a good neighbor of $S_1$ or a non-neighbor of $S_1$.

\item[$(5)$]
 if $|S_0|\ne 2$, then
$\l S_0,y_1,y_0,S_1\r $  contains no multi-multiple edges;\\
 if $|S_0|=2$, then the edge of $S_0$ has the maximum possible
multiplicity in $\l S_0,y_1,y_0,S_1\r $;

\item[$(6)$]
If  $\l y_0,S_1\r $ is either a Lann\'er diagram or a diagram of a simplicial prism 
without tail, then exactly one of the following holds:\\
$\bullet$ \ either $S_1$ is a maximal connected elliptic subdiagram in  $\l y_0,S_1\r $ of order $d-|S_0|$, \\
\smallskip
$\bullet$ \ or $S_1$ is a next to maximal connected elliptic subdiagram in  $\l y_0,S_1\r $ of order $d-|S_0|$,
$S_1$ contains a node $x$ which is an end of the dotted edge,
and the diagram $\l y_0,S_1\r \setminus x$ is a unique 
 maximal connected elliptic subdiagram of order $d-|S_0|$ in $\l y_0,S_1\r$.

\end{itemize}

\end{lemma}


%

\begin{table}[!h]
\begin{center}
\caption{List of diagrams   $\l y_0,S_1\r $, see Lemma~\ref{l1}.}
\label{list}
\vspace{10pt}
\psfrag{y}{\scriptsize $y$}
\psfrag{8}{\scriptsize $8$}
\psfrag{10}{\scriptsize $10$}
\epsfig{file=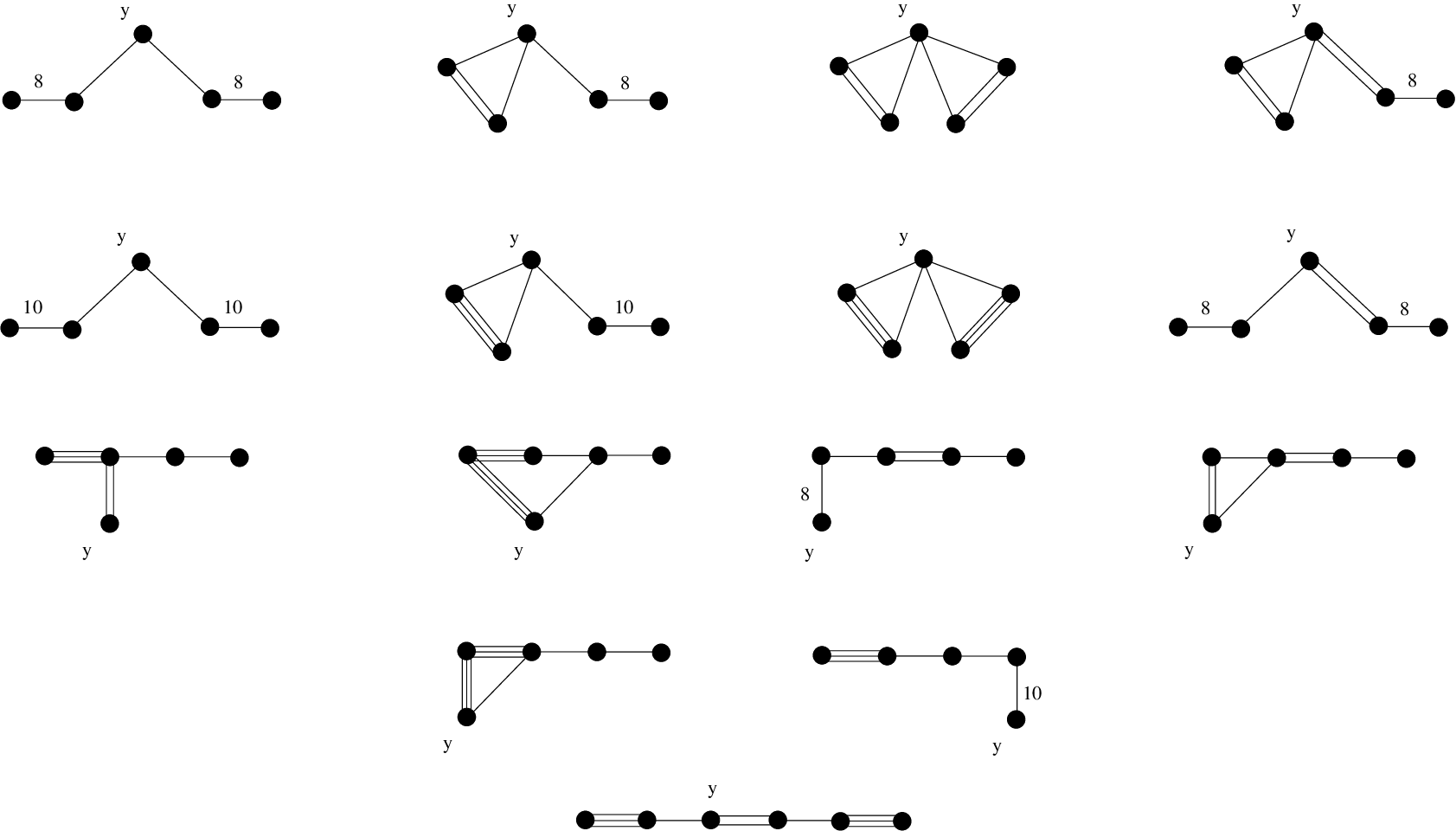,width=0.99\linewidth}
\end{center}
\end{table}

\begin{proof}

The construction is very close to one provided in~\cite[Lemma~3]{nodots}.

\begin{itemize}
\item[1.] 
{\it Analyzing the data.}
Since $S_0$ has no good neighbors,  $\overline S_0=\Sigma_{S_0}$.
Let $d_0=d-|S_0|$ be the dimension of $P(S_0)$.
Being a sub\-dia\-gram of $\Sigma(P)$,
the diagram $\Sigma_{S_0}$ contains at most one dotted edge.
Clearly, $d_0<d$. By the assumption,
 the Main Theorem holds for polytopes of dimension less than $d$,
so $P(S_0)$ contains at most $d_0+3$ facets,
and it is either a simplex, or a $d_0$-prism, or an Esselmann polytope, or a $d_0$-polytope with $d_0+3$ facets.

\item[2.] 
{\it Choosing a diagram $\Sigma'=\l S_1,y_0\r$.}

\noindent
If $P(S_0)$ is a simplex then $\Sigma'=\o S_0$.

\noindent
If $P(S_0)$ is a prism then  $\Sigma'$ is a diagram of a prism without tail.

\noindent
If $P(S_0)$ is  a $d_0$-polytope with $d_0+3$ facets then $\Sigma'$ is one of the diagrams shown in 
the first two lines of Table~\ref{list}. 

\noindent
If $P(S_0)$ is  an Esselmann polytope, then each node of $\o S_0$ belongs to some subdiagram
of the type shown in the third and fourth lines of Table~\ref{list}. Thus, we may choose as  $\Sigma'$ a
diagram of the type shown in Table~\ref{list} not containing any end of the dotted edge
(clearly, at least one such node does exist).

\item[3.] 
{\it Choosing $S_1$ and $ y_0$ in $\Sigma'$.}

\noindent
If $P(S_0)$ is an Esselmann polytope or a $d_0$-polytope with $d_0+3$ facets,
then $y_0$ is the marked node of the diagram (see Table~\ref{list}), and $S_1=\Sigma'\setminus y_0$.

\noindent
If  $P(S_0)$ is a prism, then $\Sigma'$ contains at least one connected elliptic subdiagram of order $d_0$,
and we take as $S_1$ any maximal one.

\noindent
Now, let $P(S_0)$ be a simplex. 
Consider a maximal elliptic connected subdiagram $S\subset \Sigma'$ of order $d_0$.
Let $x\subset \o S$ be a node not joined with $S_0$ by the dotted edge
(there exists one since $\o S$ is either a Lann\'er diagram or a diagram containing at least two nodes
besides $S_0$). By the choice of $x$, $\Sigma'\setminus S$ is the only node of the subdiagram $\l S_0,x,\Sigma'\r$ 
that can be joined with $x$ by the dotted edge.
If $x$ is not joined with $\Sigma'\setminus S$ by the dotted edge, 
we choose $S_1=S$ and $y_0=\Sigma'\setminus S$, otherwise we take as $S_1$  
a next to maximal  elliptic connected subdiagram of $\Sigma'$ of order $d_0$
(and $y_0=\Sigma'\setminus S_1$).

\item[4.] 
{\it Choosing $y_1$.}
Consider a subdiagram $\o S_1$. We claim that it is always possible to take a node $y_1\subset \o S_1\setminus S_0$ 
such that  $y_1$
is not joined by the dotted edge neither with $\l S_1,y_0\r$ nor with $S_0$. Indeed, we may choose $y_1$
not to be joined by the dotted edge with $S_0$ (the argument repeats one given in the preceding item).
Furthermore, such $y_1$ is not joined with  $\l S_1,y_0\r$ by the dotted edge 
 by the choice of $S_1$ and $y_0$, (see items~2 and~3).

\end{itemize}

Clearly, all conditions $(0)$--$(6)$ are satisfied by the construction.

\end{proof}

A nice property of the construction is that any edge of the obtained diagram $\l S_0,y_1,y_0,S_1\r $
belongs to either  $\l S_0,y_1\r $ or  $\l y_1,y_0,S_1\r $.
This implies that we may use the following equation on local determinants (see~\cite[Proposition~12]{abs} or~\cite[Proposition~3.1.1]{nodots}):
$$ \det(\l S_0,y_1,y_0,S_1\r ,y_1)=
\det(\l S_0,y_1\r ,y_1)+\det(\l y_1,y_0,S_1\r ,y_1)-1.$$
Combining this with the fact that $|\l S_0,y_1,y_0,S_1\r|=d+2$ 
(and, hence, $\det\l S_0,y_1,y_0,S_1\r=0 $), we obtain 
$$ \det(\l S_0,y_1\r ,y_1)+\det(\l y_1,y_0,S_1\r ,y_1)=1.
\eqno (*)
$$
This allows us to prove the finiteness of the number of diagrams  $\l S_0,y_1,y_0,S_1\r $
in consideration.

\begin{lemma}
\label{pr1fin}
The number of diagrams $\l S_0,y_1,y_0,S_1\r $
of signature $(d,1)$, $4\le d \le 8$, satisfying conditions
$(0)-(6)$ of Lemma~\ref{l1}, is finite.

\end{lemma}

\begin{proof}

It is easy to see that the number of the diagrams $\l S_0,y_1,y_0,S_1\r $ with 
$S_0\ne G_2^{(k)}$ for $k\ge 6$ is finite.
Indeed, by conditions $(0)$ and $(5)$ the diagram $\l S_0,y_1,y_0,S_1\r $
contains neither dotted nor  multi-multiple edges.
Since $|S_0|+|S_1|=d\le 8$, we obtain that
$|\l S_0,y_1,y_0,S_1\r |\le 10$,
and we have finitely many possibilities for the diagram.

Now, consider the case $S_0= G_2^{(k)}$, $k\ge 6$.
As it was mentioned above, by construction of the diagram 
 $\l S_0,y_1,y_0,S_1\r $ we may use the equation $(*)$ on local determinants.
Since $|\l y_1,y_0,S_1\r|=d$, we have
$$
|\det\l y_1,y_0,S_1\r|  \le d!
\eqno (**)
$$
(since the absolute value of each of the summands in the standard
expansion of the determinant does not exceed $1$).
Further, if $\l y_0,S_1\r $ is  a Lann\'er diagram of order 3 
then $\det\l y_0,S_1\r $ is bounded from above by $\frac{3}{4}-\cos^2(\frac{\pi}{7})\approx -0.329$
(which is the determinant of the Lann\'er diagram of order 3 with one simple edge, one empty edge
and one edge labeled by 7). 
If $\l y_0,S_1\r $ is  a diagram of a 3-prism without tail,
then the determinant of $\l y_0,S_1\r$ is a decreasing function on multiplicities 
of all edges of $\l y_0,S_1\r$. So, it is easy to check that 
 $\det\l y_0,S_1\r $ is bounded from above by $\frac{1-\sqrt 5}{16}\approx -0.08$.
In all other cases, i.e. 
if $\l y_0,S_1\r $  is neither a Lann\'er diagram of order 3 
nor a 3-prism without tail,
according to condition (3) we have finitely many possibilities for  
 $\l y_0,S_1\r $. 
Therefore, there exists a positive constant $M$ such that
$$
M<|\det\l y_0,S_1\r|.
\eqno (***)
$$
Combining $(**)$ and $(***)$, we obtain
that the determinant $\det(\!\l y_1,y_0,S_1\r ,y_1)$ (which is positive) is bounded from above,
so $\det(\l S_0,y_1\r ,y_1)$ (which is negative) is uniformly bounded from below.
However, since $S_0= G_2^{(k)}$, $k\ge 6$, the determinant $\det(\l S_0,y_1 \r,y_1)$
tends to infinity while $k$ increases (see~\cite{abs}).   
Thus, $k$ is bounded, and there are finitely many possibilities for the whole diagram 
$\l S_0,y_1,y_0,S_1\r $.

\end{proof}

Now we define several lists of diagrams similar to ones defined in~\cite[Section~3]{nodots}.

For each $S_0=G_2^{(k)},B_3,B_4,H_3,H_4,F_4$ we can write down the
complete list $$L_1(S_0,d)$$ of diagrams  $\l S_0,y_1,y_0,S_1\r $ of
signature $(d,1)$, $4\le d \le 8$, satisfying conditions $(0)-(6)$
of Lemma~\ref{l1}.
Define also a list
\begin{center}
$L_1(d)=\bigcup\limits_{k=6}^{\infty} L_1(G_2^{(k)},d).$
\end{center}
By Lemma~\ref{pr1fin}, the list $L_1(d)$ is also finite.

Clearly, the list $L_1(S_0,d)$ contains the list  $L(S_0,d)$
defined in~\cite[Section~3.2]{nodots}. 

These lists were obtained by a computer. The procedure is provided by the proof of Lemma~\ref{pr1fin}. Namely, to get the list $L(S_0,d)$ we do the following.

We list all possible diagrams $\l y_0,S_1\r $ of signature $(d-|S_0|,1)$ according to condition $(3)$ taking into account that muliplicity of an edge in $\l y_0,S_1\r $ does not exceed either $3$ (if $|S_0|\ne 2$) or $k-2$ (if $S_0=G_2^{(k)}$). For each of these diagrams we compose all possible diagrams $\l S_0,y_1,y_0,S_1\r $ by joining a new node $y_1$ with $S_0$ and $\l y_0,S_1\r $ in all possible ways by edges of multiplicity not exceeding either $3$ or $k-2$ depending on $S_0$. The list $L(S_0,d)$ consists of those diagrams $\l S_0,y_1,y_0,S_1\r $ which have zero determinant and contain no parabolic subdiagrams. 

To get the list  $L_1(d)$ we take the union of the lists $ L_1(G_2^{(k)},d)$ for $6\le k\le k_0$, where $k_0$ can be found according to the proof of Lemma~\ref{pr1fin}. More precisely, the expression for $\det(\l G_2^{(k)},y_1\r ,y_1)$ (see e.g.~\cite[Section~3.1]{nodots}) shows that for $k\ge 7$ the local determinant $\det(\l G_2^{(k)},y_1\r ,y_1)$ does not exceed $1-1/(4\sin^2{\frac{\pi}{k}})$. Combining inequalities $(**)$ and $(***)$, we see that the local determinant $\det(\!\l y_1,y_0,S_1\r ,y_1)$ is bounded from above by some constant $d!/M$ depending on $d$ only. Now, combining this with equation $(*)$, we get an explicit expression for $k_0$.

Usually the lists $L_1(S_0,d)$ and $L_1(d)$ are not very short. In what follows we reproduce
some parts of the lists as far as we need.

In fact, the bounds obtained in the proof of Lemma~\ref{pr1fin}
are not optimal. To reduce computations  we usually analyze concrete data.
For example, instead of taking $d!$ as the bound of $|\det\l S_0,y_1,y_0,S_1\r|$,
we may bound it by the number of negative summands in its expansion.
This leads to reasonable bounds on the multiplicity of multi-multiple edges in $\l S_0,y_1,y_0,S_1\r $,
the worst of which was $87$ in one of the cases.

Now, given a diagram $\Sigma$, a constant $C$ and dimension $d$,
define a list $$L'(\Sigma,C,d)$$ of diagrams $\l \Sigma,x\r$ of signature $(d,1)$
containing no subdiagrams of  the type $G_2^{(k)}$ for $k>C$,
no dotted edges incident to $x$, and no parabolic subdiagrams.
Clearly, for given $\Sigma$, $C$ and $d$, this list is finite.
One can notice that 
if $\Sigma$ contains no dotted edges, this list coincides with 
the list  $L'(\Sigma,C,d)$ defined in~\cite[Section~3.2]{nodots}.

The list $L'(\Sigma,C,d)$ is easy to obtain by computer. We join a new node with all nodes of $\Sigma$ by edges of multiplicity at most $C-2$ and choose those diagrams having signature $(d,1)$ and containing no parabolic subdiagrams. To reduce the computations, we first compute the determinant, and check the signature only if the determinant vanishes. 

As in~\cite{nodots}, for  given $\Sigma$, $C$, $d$ and an elliptic subdiagram $S\subset \Sigma$
we define the sublist   $L'(\Sigma,C,d,S^{})$ which consists of diagrams  $\l \Sigma,x\r$,
where either $x$ is either a good neighbor or a non-neighbor of $S$ (in~\cite{nodots} this list 
is denoted by $L'(\Sigma,C,d,S^{(g,n)})$.

\section{Proof of the Main theorem} 
\label{proof}

First, we prove some general facts concerning subdiagrams of the type $B_k$ which will be
used later for the proof in all dimensions;
then we prove the theorem starting from low dimensions and going to higher ones.

\subsection{Subdiagrams of the type $B_k$}

\begin{lemma}
\label{b_d}
Let $P\subset \H^d$, $d\ge 6$, be a compact Coxeter polytope such that $\Sigma(P)$ contains a unique dotted edge.
If $\Sigma(P)$ 
contains neither subdiagram of the type $F_4$ nor subdiagram of the type
$G_2^{(k)}$, $k\ge 5$,
then   $\Sigma(P)$ contains no subdiagram of the type $B_d$.

\end{lemma}

\begin{proof}
At first, notice that the assumptions of the lemma imply that for any two nodes $t_i,t_j\in \Sigma$
we have $[t_i,t_j]\in \{2,3,4,\infty\}$ (recall that $[t_i,t_j]=k$ means that the nodes $t_i$ and $t_j$ are joined by a $(k-2)$-fold edge, and $[t_i,t_j]=\infty$ means that the nodes are joined by a dotted edge). This will be used frequently throughout the proof. 

Suppose that $S_0\subset  \Sigma(P)$ is a diagram of the type $B_d$,
denote by $t_1,\dots,t_d$ the nodes of $S_0$ ($[t_1,t_2]=4$, $[t_i,t_{i+1}]=3$ for all $1<i<d$).
Consider the diagram $S_1=\l t_1,t_2,\dots,t_{d-1}\r$ of the type $B_{d-1}$.
The polytope $P(S_1)$ is one-dimensional, so the diagram $\Sigma_{S_1}$ consists of two nodes connected by a dotted edge.
By~\cite[Theorem~2.2]{All}, this implies that the diagram $\l S_1,\o S_1\r $ is of one of the two types 
shown in Fig.~\ref{Bd_ab} (since $t_d\in\o S_1$). We consider these two diagrams separately.

\begin{figure}[!h]
\begin{center}
\psfrag{a}{(a)}
\psfrag{b}{(b)}
\psfrag{x}{\small $x$}
\psfrag{1}{\small $t_1$}
\psfrag{2}{\small $t_2$}
\psfrag{3}{\small $t_3$}
\psfrag{4}{\small $t_4$}
\psfrag{d2}{\small $t_{d-2}$}
\psfrag{d1}{\small $t_{d-1}$}
\psfrag{d}{\small $t_{d}$}
\psfrag{2,3}{\scriptsize $2,3$}
\psfrag{dots}{$\dots$}
\epsfig{file=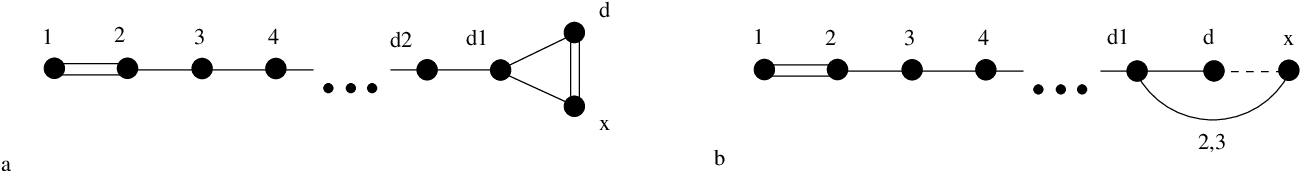,width=0.95\linewidth}
\caption{Two types of the diagram  $\l S_1,\o S_1\r $, see Lemma~\ref{b_d}}
\label{Bd_ab}
\end{center}
\end{figure}

\medskip
\noindent
{\bf Case (1):} $\l S_1,\o S_1\r $ is a diagram of the type 
shown in Fig.~\ref{Bd_ab}(a).\\ 
Consider the diagram $S_2=\l t_2,t_3,\dots,t_{d-1},t_d \r$
of the type $A_{d-1}$. It has a unique good neighbor in  $\l S_1,\o S_1\r $,
so in $\Sigma$ there exists a node $y$ which is either a good neighbor or a non-neighbor of $S_2$
(since the diagram of the type $A_{d-1}$ defines a 1-face of $P$, which should have two ends).
We consider two cases: either $y$ is joined with $t_1$ by a dotted edge, or it is not.

\smallskip
\noindent
{\bf Case (1a):} Suppose that $[y,t_1]=\infty$.\\ 
Consider the diagram $S_3=\l t_1,t_2,\dots,t_{d-3}\r$ of the type $B_{d-3}$. $P(S_3)$ is a Coxeter 3-polytope
whose Coxeter diagram $\Sigma_{S_3}$ contains a Lann\'er subdiagram of order 3 (coming from the subdiagram 
$\l t_{d-1}, t_d,x \r\subset  \Sigma$). This implies that  $P(S_3)$ is not a simplex, so, it has 
a pair of non-intersecting facets. Since  $\Sigma$ contains only one dotted edge $yt_1$, which is not contained in $\o S_3$, 
we conclude that $S_3$ has a good neighbor $z$. So, $z$ is not joined with $\l t_1,t_2,\dots,t_{d-4}\r$,
$[z,t_{d-3}]=3$ (here we use that $d\ge 6$ and that $\Sigma$ contains no subdiagram of the type $F_4$).
Furthermore, $z$ may be joined with $t_{d-1}$, $t_d$ and $x$ by either simple or double edge.  
Notice, that
$[z,t_{d-2}]=4$, otherwise either $\l t_{d-3},t_{d-2},z \r$ or 
$\l S_3,t_{d-2},z \r$  
is a parabolic subdiagram (of the types $\w A_2$ and $\w B_{d-2}$ respectively).
So, we have 27 possibilities for the diagram $\l S_0,x,z \r=\l t_1,t_2,\dots,t_{d-1},t_d,x,z\r$
(see Fig.~\ref{Bd_a}(a)). The diagram $\l S_0,x,z \r$ contains $d+2$ nodes, so 
$\det \l S_0,x,z \r=0$, which holds only in the case shown in  Fig.~\ref{Bd_a}(b) 
(to see this for any $d\ge 5$, we use local determinants, namely, we check the equality 
$\det(S_3,t_{d-3})+\det(\l x,z,t_{d-3},t_{d-2},t_{d-1},t_d \r,t_{d-3})=1$).  

Recall that $y$ is either a good neighbor or a non-neighbor of $S_2$.
So, $y$ is joined with $S_2$ by at most one edge (simple or double, since $[y,t_1]=\infty$).
On the other hand, $y$ should be joined with each of the Lann\'er diagrams
$\l z,t_{d-3},t_{d-2}\r$, $\l z,t_{d-2},t_{d-1}\r$ and $\l x,t_{d-1},t_d\r$.
Since any non-dotted edge in $\Sigma(P)$ has multiplicity at most two,
a short explicit check shows that we always obtain a parabolic subdiagram of one of the types $\w A_2$, $\w C_2$, 
$\w A_3$ and $\w C_3$, which is impossible.

\begin{figure}[!h]
\begin{center}
\psfrag{a}{(a)}
\psfrag{b}{(b)}
\psfrag{c}{(c)}
\psfrag{x}{\scriptsize $x$}
\psfrag{y}{\scriptsize $y$}
\psfrag{z}{\scriptsize $z$}
\psfrag{1}{\scriptsize     $t_1$}
\psfrag{2}{\scriptsize     $t_2$}
\psfrag{3}{\scriptsize     $t_3$}
\psfrag{4}{\scriptsize     $t_4$}
\psfrag{d4}{\scriptsize     $t_{d-4}$}
\psfrag{d3}{\scriptsize     $t_{d-3}$}
\psfrag{d2}{\scriptsize     $t_{d-2}$}
\psfrag{d1}{\scriptsize     $t_{d-1}$}
\psfrag{d}{\scriptsize     $t_{d}$}
\psfrag{2,3,4}{\tiny $2,\!3,\!4$}
\psfrag{3,4}{\tiny $3,\!4$}
\psfrag{dots}{$\dots$}
\epsfig{file=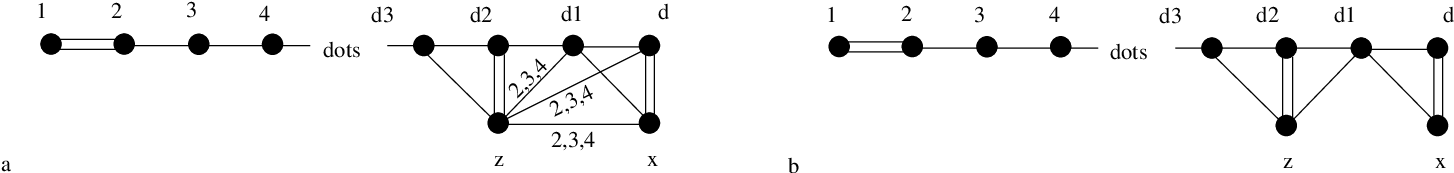,width=0.97\linewidth}
\caption{To the proof of Lemma~\ref{b_d}, case (1a). }
\label{Bd_a}
\end{center}
\end{figure}

\smallskip
\noindent
{\bf Case (1b):} Suppose that $[y,t_1]\ne \infty$. \\
Since $y$ is either a good neighbor or a non-neighbor of $S_2=\l t_2,\dots,t_{d-1},t_d \r$,
$y$ cannot be joined with $S_0$ by a dotted edge. However, it is possible that $[y,x]=\infty$.
In the latter case we consider the diagram $S_2'=\l t_2,\dots ,t_{d-1},x\r$
instead of the diagram $S_2$ and find its good neighbor (or non-neighbor) $y'$,
which is definitely not an end of a dotted edge in this case.
Therefore, we may assume that $[y,x]\ne \infty$, in other words,   
that the diagram $\l S_0,x,y\r$ contains no dotted edges.

To find out, how $y$ can be joined with $\l S_0,x\r$, notice that:
\begin{itemize}

\item[1.]
$y$ is joined with $S_0$ and with $\l S_1,x\r$ (otherwise we obtain an elliptic diagram of order $d+1$).
\item[2.]
$[y,t_1]\ne 2$ (otherwise either the subdiagram $\l S_0,y\r$ contains a parabolic subdiagram, 
or $\l S_0,y\r$ is a diagram of the type $B_{d+1}$, which is also impossible). 
\item[3.] 
$y$ is joined with one of $t_2$
and $t_3$ (otherwise $\l y,t_1,t_2,t_3\r$ either is a diagram of the type $F_4$ or contains a parabolic subdiagram
of the type $\w C_2$). In particular, this implies that $y$ is not joined with $\l t_4,t_5,\dots,t_d\r$.
\item[4.]  
$[y,x]\ne 2$ (since  $d\ge 6$, the edge $yx$ is the only way to join an indefinite subdiagram $\l y,t_1,t_2,t_3\r$ 
with Lann\'er diagram $\l t_{d-1},t_d,x \r$).
\item[5.] 
$[y,x]=3$
(if $[y,x]=4$
then $\l y,x,t_d\r$ is a parabolic diagram of the type $\w C_2$).
\item[6.] 
$[y,t_1]=3$
(if $[y,t_1]=4$
then $\l t_1,y,x,t_d\r$ is a parabolic diagram of the type $\w C_3$).
\item[7.] 
$[y,t_2]=3$
(if $[y,t_2]=2$ then $\l t_2,t_1,y,x,t_d\r$ is a parabolic diagram of the type $\w C_4$,
if $[y,t_2]=4$ then $\l t_2,y,x,t_d\r$ is a parabolic diagram of the type $\w C_3$.   

\end{itemize}

We arrive with a parabolic subdiagram $\l x,y,t_2,t_3,t_4,\dots,t_{d-2},t_{d-1} \r$
of the type $\w A_{d-1}$, which is impossible.

\medskip
\noindent
{\bf Case (2):} $\l S_1,\o S_1\r $ is a diagram of the type 
shown in Fig.~\ref{Bd_ab}(b).\\
Similarly to the case (1), we consider the diagrams $S_2=\l t_2,t_3,\dots,t_{d-1},t_d \r$
and $S_3=\l t_1,t_2,\dots,t_{d-2}\r$.
As before, $S_2$ has either a good neighbor or a non-neighbor $y$,
and $S_3$ has a good neighbor $z$ (to see the latter statement, notice, that $P(S_3)$ is a 2-polytope 
whose diagram $\Sigma_{S_3}$ contains a dotted edge coming from $\l t_d,x\r$,
so  $\Sigma_{S_3}$ contains at least one more dotted edge, which can appear only if $S_3$
has one more good neighbor). So, $[z,t_{d-2}]=3$, which implies   $[z,t_{d-1}]=4$
(otherwise, either $\l S_3,t_{d-1},z\r$ is a parabolic diagram of the type $\w B_{d-1}$ or $\w C_{d-1}$, 
or $\l t_{d-2},t_{d-1},z\r$ is of the type $\w A_{2}$).
So, $\l S_0,z\r$ is one of the two diagrams shown in  Fig.~\ref{bd_b}(a).

\begin{figure}[!h]
\begin{center}
\psfrag{a}{(a)}
\psfrag{b}{(b)}
\psfrag{c}{(c)}
\psfrag{x}{\scriptsize $x$}
\psfrag{y}{\scriptsize $y$}
\psfrag{z}{\scriptsize $z$}
\psfrag{1}{\scriptsize $t_1$}
\psfrag{2}{\scriptsize $t_2$}
\psfrag{3}{\scriptsize $t_3$}
\psfrag{4}{\scriptsize $t_4$}
\psfrag{5}{\scriptsize $t_5$}
\psfrag{6}{\scriptsize $t_6$}
\psfrag{d4}{\scriptsize $t_{d-4}$}
\psfrag{d3}{\scriptsize $t_{d-3}$}
\psfrag{d2}{\scriptsize $t_{d-2}$}
\psfrag{d1}{\scriptsize $t_{d-1}$}
\psfrag{d}{\scriptsize $t_{d}$}
\psfrag{2,3,4}{\tiny $2,3,4$}
\psfrag{2,3}{\tiny $2,3$}
\psfrag{2.3}{\tiny $2,\!3$}
\epsfig{file=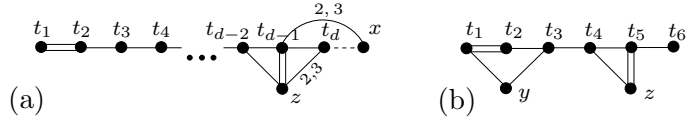,width=0.7\linewidth}
\caption{To the proof of Lemma~\ref{b_d}, case (2).}
\label{bd_b}
\end{center}
\end{figure}

Similarly to case (1b), consider the multiplicities of edges joining $y$ with  $\l S_1,z\r$. 
All the assertions 1--7 (as well as the arguments) still hold 
if we replace $x$ by $z$, $t_{d}$ by $t_{d-1}$, and  $t_{d-1}$ by $t_{d-2}$.
However, to state assertion 4 we need to assume now that $d\ge 7$.
To state the same for $d=6$ notice, that the only case when $[y,z]=2$ and all Lann\'er subdiagrams of
$\l y,t_1,t_2,t_3\r$ are joined
with Lann\'er diagram $\l t_{d-2},t_{d-1},z \r$ is one shown in Fig.~\ref{bd_b}(b)
 (in all other cases the subdiagram 
$\l t_1,\dots,t_5,y,z\r$ contains a parabolic subdiagram).
However, this diagram is superhyperbolic, so all the assertions 1--7 hold for any $d\ge 6$.
This leads to a parabolic subdiagram $\l z,y,t_2,t_3,t_4,\dots,t_{d-3},t_{d-2} \r$
of the type $\w A_{d-2}$, which is impossible. 

\end{proof}

\begin{lemma}
\label{b_k}
Let $P\subset \H^d$, $d\ge 4$, be a compact Coxeter polytope such that $\Sigma(P)$ contains a unique dotted edge.
Suppose that $\Sigma(P)$ 
contains no subdiagram of the type $F_4$, 
$G_2^{(m)}$, $m\ge 5$, and  $B_d$.
Then
$\Sigma(P)$ contains no subdiagram of the type $B_k$ for any $k< d$, $k\ge 3$.

\end{lemma}

\begin{proof}
Suppose that the lemma is true for all $k'>k$, but there exists a subdiagram
$S_0\subset \Sigma(P)$ of the type $B_k$. 
Then $S_0$ has no good neighbors (here we use the assumption that $\Sigma$ contains no subdiagram of the type $F_4$).
Thus, $\o S_0=\Sigma_{S_0}$ is a Coxeter diagram of a $(d-k)$-polytope $P(S_0)$.
Clearly, $\o S_0$ contains at most one dotted edge and does not contain edges of multiplicity greater 
than 2. 
As above, denote by $t_1,t_2,\dots,t_d$ the nodes of $S_0$ ($[t_1,t_2]=4$, $[t_i,t_{i+1}]=3$ for all $1<i<d$).

Consider a subdiagram $S_1\subset S_0$ of the type $B_{k-1}$. Since $S_1\subset S_0$, at least one bad neighbor 
(denote it by  $x$) of $S_0$
is not a bad neighbor of $S_1$ ($P(S_1)$ is a face of bigger dimension than $P(S_0)$ is).
Suppose that $x$ is not an end of the dotted edge.  
Clearly, $x$ is a good neighbor of $S_1$, otherwise it is a non-neighbor and the diagram 
$\l S_0,x\r$ is either a parabolic diagram $\w C_{k}$ or a diagram of the type $B_{k+1}$ which is impossible 
by assumption. So, $\l S_1,x\r$ is a diagram of the type $B_k$
(we use the assumption that $k> 3$ and that $\Sigma(P)$ contains no subdiagram of the type $F_4$).
Let $x'$ be any node of $\o S_0$ joined with $x$ (it does exist since an indefinite diagram 
$\l S_0,x\r$ should be joined with each Lann\'er subdiagram of $\o S_0$).
Then the diagram $\l S_1,x,x'\r$ is either a parabolic diagram $\w C_{k}$ or a diagram of the type 
$B_{k+1}$, which is impossible by assumption. 

Therefore, $x$ is an end of the dotted edge.
Moreover, the paragraph above shows that another  end of the dotted edge
coincides with either $t_d$ or some $x'\subset \o S_0$ (otherwise we repeat the arguments and obtain a contradiction).
This implies that $x$ is the only bad neighbor of $S_0$ that is not a bad neighbor of $S_1$, 
and either $[x,t_k]=\infty$ or $[x,x']=\infty$, where $x'\in \o S_0$. 
In particular, this implies that $\o S_0$ contains no dotted edge, which is possible only if
$\Sigma_{S_0}$ is one of the diagrams shown in Fig.~\ref{bk}
(here we use the classification of Coxeter polytopes with mutually intersecting facets,
we also use that any non-dotted edge of $\Sigma$ is either a simple edge or a double edge). 

\begin{figure}[!h]
\begin{center}
\epsfig{file=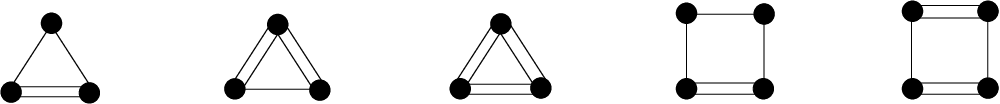,width=0.7\linewidth}
\caption{Possible diagrams $\Sigma_{S_0}=\o S_0$, see Lemma~\ref{b_k}.}
\label{bk}
\end{center}
\end{figure}

Suppose that $[x,x']=\infty$, where $x'\in \o S_0$. It is easy to see that 
 $[x,t_{k-1}]=3$ and  $[x,t_k]=4$ (otherwise $\Sigma$ contains either a parabolic subdiagram or a
subdiagram of the type $B_{k+1}$). Since $x$ is the only bad neighbor of $S_0$ that is not a bad neighbor of $S_1$,
we have $\o S_1=\l t_k,x, \o S_0 \r$. Thus, the diagram $\Sigma_{S_1}$ 
contains exactly three Lann\'er subdiagrams:
two dotted edges coming from $t_kx$ and $xx'$, and a Lann\'er diagram of order 2 or 3 (which coincides with $\o S_0$).
Hence, the Lann\'er diagram coming from $xx'$ has a common point with any other Lann\'er diagram
of  $\Sigma_{S_1}$, which is impossible by~\cite[Lemma~1.2]{d+4}.

Therefore, $[x,t_k]=\infty$.
Let $S_2=\l t_2,t_3,\dots,t_k \r$ be a subdiagram of the type $A_{k-1}$,
and let $S_3\subset \o S_0$ be any subdiagram of the type $B_3$ (if any) or of the type $B_2$
(otherwise). Then the subdiagram $\l S_2,S_3\r$ 
has exactly one good neighbor  (or non-neighbor) $y$ besides the node $t_1$.
Clearly, $y$ is a bad neighbor of $S_0$ distinct from $x$. So, $y$ is not an end of the dotted edge.
Let $y'=\o S_0\setminus S_3$.

\noindent
To find out, how $y$ can be joined with $\l S_0,x\r$, notice that:
\begin{itemize}
\item[1.]
$[y,t_1]\ne 2$ (otherwise the subdiagram $\l S_0,y\r$ contains a parabolic subdiagram). 
\item[2.] 
$y$ is joined with one of $t_2$
and $t_3$ (otherwise $\l y,t_1,t_2,t_3\r$ either is a diagram of the type $F_4$, 
or contains a parabolic subdiagram
of the type $\w C_2$). In particular, this implies that $y$ is not joined with $\l t_4,t_5,\dots,t_k\r$.
\item[3.]  
$y$ is not joined with $S_3$ (otherwise an elliptic diagram $\l S_2,S_3,y\r$ is connected, so it is 
of the type $B_{k+2}$ or $B_{k+3}$).
\item[4.]
$[y,y']=3$
(if $[y,y']=4$
then $\l t_1,y,\o S_0\r$ contains either a parabolic diagram of the type $\w C_2$ or $\w C_3$,  
or a subdiagram of the type $F_4$).

\end{itemize}

Either $\l t_1,t_2,t_3,y\r$ or $\l t_1,t_2,y\r$ is a Lann\'er diagram 
(one of the diagrams shown in Fig.~\ref{bk}), denote it by $L$.
By construction, $L$ is joined with a Lann\'er diagram $\o S_0$  
by the edge $yy'$ only. Thus, we obtain a subdiagram $\l L,\o S_0\r\subset\l S_0,y,\o S_0\r$ of the following type:
it consists of two Lann\'er diagrams $L$ and $\o S_0$ from Fig.~\ref{bk} joined by a unique simple edge  $yy'$, 
where $y\in L$, $y'\in \o S_0$, and both diagrams $L\setminus y$ and $\o S_0\setminus y'$ are of the type $B_2$ or $B_3$.
It is easy to see that any such diagram $\l L,\o S_0\r$ is superhyperbolic,
which proves the lemma.

\end{proof}

\begin{lemma}
\label{b_2}
Suppose that the Main Theorem holds for any dimension $d'<d$, $d> 4$.
Suppose also that for any compact Coxeter polytope $P\subset \H^d$, such that $\Sigma(P)$ contains a unique dotted edge, 
it is already shown that $\Sigma(P)$ contains neither subdiagram of the type $F_4$, nor subdiagram of the type
$G_2^{(k)}$, $k\ge 5$, nor subdiagram of the type $B_d$. 
Then the Main Theorem holds in dimension $d$.

\end{lemma}

\begin{figure}[!h]
\begin{center}
\psfrag{2,3}{\scriptsize $2,\!3$}
\epsfig{file=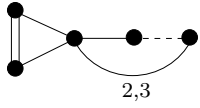,width=0.2\linewidth}
\caption{The diagram $\Sigma_{S_0}=\o S_0$, see Lemma~\ref{b_2}.}
\label{b_2_5}
\end{center}
\end{figure}

\begin{proof}
Suppose that the Main Theorem is broken in dimension $d$.
Let $P\subset \H^d$ be a compact Coxeter polytope with at least $d+4$ facets, 
such that $\Sigma(P)$ contains a unique dotted edge, 
and $\Sigma(P)$ contains neither subdiagram of the type $F_4$ nor subdiagram of the type
$G_2^{(k)}$, $k\ge 5$, nor subdiagram of the type $B_d$. 
By Lemma~\ref{b_k}, $\Sigma(P)$ also contains no subdiagram of the type $B_k$, $k>2$.
It follows that any Lann\'er diagram of $\Sigma(P)$ is either a dotted edge or 
 one of the three diagrams of order three shown in  Fig.~\ref{bk}.

Let $L_0\subset \Sigma(P)$ be a Lann\'er diagram of order 2, i.e. a dotted edge.
By~\cite[Lemma~1.2]{d+4}, $\Sigma(P)\setminus L_0$ contains at least one Lann\'er diagram $L$.
So, $L$ is  one of three diagrams of order three shown in   Fig.~\ref{bk}.
Let $S_0\subset L$ be a subdiagram of the type $B_2$.
By assumptions, $S_0$ has no good neighbors, so $\o S_0=\Sigma_{S_0}$ is a diagram containing at most one dotted edge.
$\o S_0$ is a diagram of a $(d-2)$-polytope with at most $(d-2)+3$ nodes,
containing no edges of multiplicity greater than 2, and no diagrams of type $B_3$. 
It follows from the classification of $d'$-polytopes with at most $d'+3$ facets, that
$P(S_0)$ is a polytope of dimension at most 3. If $P(S_0)$ is either a 2-polytope or an 1-polytope, 
then $d<5$ in contradiction to the assumptions. 

So, $P(S_0)$ is a 3-polytope. Then  $P(S_0)$ is a 3-prism
(it cannot be a simplex since diagrams of 3-simplices 
 always contain subdiagrams of one of the forbidden types).
It is easy to see that $\o S_0=\Sigma_{S_0}$ is the diagram shown in Fig.~\ref{b_2_5}. 
Since the 5-polytope $P$ has at least $d+4=9$ facets, there exists a node $x\in \Sigma(P)$
such that $x\notin \l S_0,\o S_0\r$.
Notice that $x$ is joined with $ \l S_0,\o S_0\r$ by simple and double edges only.
Since $P$ is a 5-polytope, $\det \l x,S_0,\o S_0\r=0$. However, each of the diagrams satisfying 
all the conditions above either contains a parabolic subdiagram, or is superhyperbolic 
(in other words, the list $L'(\l S_0,\o S_0\r,4,5)$ is empty).
This proves the lemma.

\end{proof}

%
%
%
%
%
%
%

\subsection{Dimensions 2 and 3}

In dimensions 2 and 3 the statement of the Main Theorem is combinatorial:
it is easy to see that
any polygon except triangle has at least two pairs of disjoint sides, and
any polyhedron (3-polytope) having a unique pair of disjoint facets is a triangular prism.

\subsection{Dimension 4}

Let $P$ be a 4-dimensional compact hyperbolic Coxeter polytope such that $\Sigma(P)$ contains a unique dotted edge and 
$P$ has at least 8 facets.

\begin{lemma}
\label{4_mult}
$\Sigma(P)$  contains no multi-multiple edges.

\end{lemma}

\begin{proof}
Suppose that $S_0\subset \Sigma(P)$ is a multi-multiple edge of the maximum multiplicity in $\Sigma(P)$. 
Then $S_0$ has no good neighbors and, by Lemma~\ref{l1}, $\Sigma(P)$ contains a subdiagram $\l S_0,y_1,y_0,S_1 \r$ from 
the list $L_1(4)$. The list contains 28 diagrams, 3 of these diagrams are Esselmann diagrams, which cannot be subdiagrams of
$\Sigma(P)$ by~\cite[Lemma~1]{nodots}. For each of the remaining 25 diagrams we check the list 
$L'( \Sigma_1,k(\Sigma_1),4)$, where $\Sigma_1$ ranges over the 25 diagrams,
and $k(\Sigma_1)$ is the maximum multiplicity of an edge in $\Sigma_1$ (in fact,  $k(\Sigma_1)\le 14$; 
$\Sigma(P)$ contains some diagram from one of these lists by Lemma~\ref{l'}).
All these lists are empty, so the lemma is proved.  

\end{proof}

In the proof of the following lemma we use Gale diagram of simple polytope 
(see~\cite[Section~2.2]{nodots} for essential facts about Gale diagrams, and~\cite{G} for general theory).

\begin{lemma}
\label{3plus3}
$\Sigma(P)$ contains two non-intersecting Lann\'er diagrams of order 3, all nodes of which are not ends of the dotted edge.

\end{lemma}

\begin{proof}
The proof follows the proof of~\cite[Lemma~8]{nodots}.
 
Let $n$ be the number of facets of $P$ and let 
$f_{n-1} $ and $f_n$ be the facets of $P$ having no common point.

Let $\cal G$ be a Gale diagram of $P$.
It consists of $n$ points $a_1,\dots,a_n$ in $(n-6)$-dimensional sphere $\S^{(n-6)}$.
Let $a_i$  be the point corresponding to facet $f_i$.
Consider a unique hyperplane $H\subset \S^{(n-6)}$ containing all points $a_i$, $i\ge 7$.
Let  $H^+$ and $H^-$ be open hemispheres of $\S^{(n-6)}$ bounded by $H$. 
Since any two of $f_j$, $1\le j\le 6$, have non-empty intersection, each of 
$H^+$ and $H^-$ contains at least three points $a_j$, $1\le j \le 6$.
Since $n\ge 8$, $H^+$ and $H^-$ do not contain neither  $a_{n-1}$ nor  $a_n$, which proves the lemma.

\end{proof}

\begin{lemma}
\label{4,b}
The Main Theorem holds in the dimension $d=4$. 

\end{lemma}

\begin{proof}
Suppose that the Main Theorem does not hold for  $d=4$, so let $P$ be a compact Coxeter 4-polytope
with at least $8$ facets such that   $\Sigma(P)$ contains a unique dotted edge.

By Lemma~\ref{3plus3}, $\Sigma(P)$ contains two disjoint  Lann\'er subdiagrams $T_1$ and $T_2$ of order three
each such that the diagram $\l T_1,T_2\r$ contains no dotted edges.
It is shown in~\cite[Lemma~9]{nodots} that there are only 39 diagrams  $\l T_1,T_2\r$ of signature $(4,1)$ such that  
 $T_1$ and $T_2$ are Lann\'er diagrams of order three and   $\l T_1,T_2\r$ contains no edges of multiplicity 
greater than three. 3 of these diagrams are Esselmann diagrams (by~\cite[Lemma~1]{nodots}, they are not parts of
any diagram of a 4-polytope with more than 6 facets), 5 of them contain parabolic subdiagrams.
For each of the remaining 31 diagrams the list $L'(\l T_1,T_2\r, 5,4)$ is empty.

\end{proof}

\subsection{Dimension 5}

Let $P$ be a 5-dimensional compact hyperbolic Coxeter polytope such that $\Sigma(P)$ contains a unique dotted edge and 
$P$ has at least 9 facets.

\begin{lemma}
\label{5_mult}
$\Sigma(P)$  contains no multi-multiple edges.

\end{lemma}

\begin{proof}
Suppose that $S_0\subset \Sigma(P)$ is a multi-multiple edge of the maximum multiplicity in $\Sigma(P)$. 
Then $S_0$ has no good neighbors and, by Lemma~\ref{l1}, $\Sigma(P)$ contains a subdiagram $\l S_0,y_1,y_0,S_1 \r$ from 
the list $L_1(5)$. The list consists of 11 diagrams shown in Table~\ref{5_multi}.
Notice that $\o S_0$ in this case is a diagram of a 3-polytope with at most one pair of non-intersecting facets,
i.e. either a simplex or a prism. 
In the cases when $S_1$ is either a diagram of a prism without tail or a next to maximal subdiagram of a diagram
of a simplex, we mark the end of the dotted edge by a  circle. 
Denote by $S_2$ an elliptic subdiagram of $\l S_0,y_1,y_0,S_1 \r$ of order 4 marked by a gray block
(if any, see Table~\ref{5_multi}).
Notice that $S_2$  has at most 1 good neighbor or non-neighbor in $\l S_0,y_1,y_0,S_1 \r$, 
and if it has exactly one then $S_2$ contains an end of the dotted edge. Therefore, there exists a node
 $x\in \Sigma(P)\setminus \l S_0,y_1,y_0,S_1 \r$ such that $x$ is not a bad neighbor of $S_2$, and
the diagram $\l x,S_0,y_1,y_0,S_1 \r$ contains no dotted edges. In other words,
$\Sigma(P)$ contains a diagram from the list $L'(\Sigma_1,k(\Sigma_1),5,S_2)$, where
$\Sigma_1$ ranges over the 11 diagrams $\l S_0,y_1,y_0,S_1 \r$ and $k(\Sigma_1)$ 
is a maximum multiplicity of the edge in $\Sigma_1$ (in a unique case when the diagram $S_2$ is not defined, 
we take a list $L'(\Sigma_1,10,5)$ instead).
All these lists but one are empty. The remaining one contains a unique entry $\Sigma_2$ shown in  Fig.~\ref{5_multi2} 
(again, we mark an end of the dotted edge by a circle).
Consider a subdiagram $S_3\subset \Sigma _2$ of the type $G_2^{(8)}$ marked in Fig.~\ref{5_multi2} by a gray block.
Clearly, the subdiagram $\o S_3$ contains no dotted edges. At the same time, starting from $S_3$ instead of $S_0$, 
we should obtain some diagram of the list  $L_1(S_3,5)\subset L_1(5)$, but looking at Table~\ref{5_multi} one can note that 
each entry of  $L_1(5)$ containing the subdiagram $G_2^{(8)}$ contains an end of the dotted edge. 
The contradiction proves the lemma.

\begin{table}[!t]
\begin{center}
\caption{The list $L_1(5)$. Ends of dotted edges are encircled.}
\label{5_multi}
\vspace{10pt}
\psfrag{a}{(a)}
\psfrag{b}{(b)}
\psfrag{6}{\scriptsize $6$}
\psfrag{8}{\scriptsize $8$}
\psfrag{10}{\scriptsize $10$}
\psfrag{12}{\scriptsize $12$}
\epsfig{file=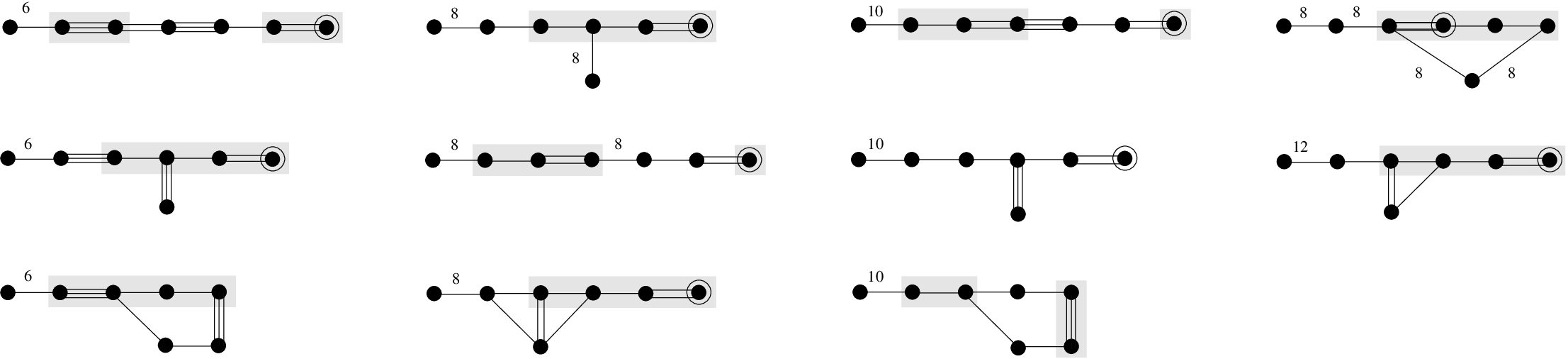,width=0.997\linewidth}
\end{center}
\end{table}

\begin{figure}[!h]
\begin{center}
\psfrag{a}{(a)}
\psfrag{b}{(b)}
\psfrag{6}{\small $6$}
\psfrag{8}{\scriptsize $8$}
\psfrag{10}{\small $10$}
\psfrag{12}{\small $12$}
\epsfig{file=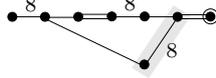,width=0.22\linewidth}
\caption{Treating the list $L_1(5)$, see Lemma~\ref{5_mult}. }
\label{5_multi2}
\end{center}
\end{figure}

\end{proof}

\begin{lemma}
\label{5_h4}
$\Sigma(P)$  contains no subdiagrams of the types $H_4$.

\end{lemma}

\begin{proof}
Suppose that $S_0\subset \Sigma(P)$ is a subdiagram of the type $H_4$.
Then $S_0$ has no good neighbors, so $\o S_0=\Sigma_{S_0}$ is a dotted edge.
Let $S_1\subset S_0$ be a subdiagram of the type $H_3$.
By Lemma~\ref{polygon}, $S_1$ has a good neighbor or a non-neighbor $x\notin \l S_0,\o S_o \r$.
If $x$ is a good neighbor of $S_1$,
consider the diagram $S_2=\l S_1,x\r$ of the type $H_4$. As it is shown above for the 
diagram $S_0$, the dotted edge belongs to $\o S_2$. 
Hence, the dotted edge is not joined with an indefinite diagram $\l S_0,x\r$, which is impossible.  
Therefore, $x$ is a non-neighbor of $S_1$. Let $y$ be an end of the dotted edge joined with $x$
(there exists one, since $\Sigma(P)$ is not superhyperbolic).
Let $t_1=S_0\setminus S_1$ and notice that $[x,t_1]\ne 5$ (otherwise $\l S_0,x \r$ contains a subdiagram $S_3$ of the type $H_4$
such that $\o S_3$ contains no dotted edge, which is impossible as it was proved above).
Thus, we have only 6 possibilities for the diagram $\l S_0,x,y\r$ (see Fig.\ref{5_h4_}(a)). 
In fact, only in 3 of these cases the diagram  $\l S_0,x,y\r$ contains no parabolic subdiagrams.
If $x$ is joined with $S_0$ by a simple edge, we consider the list $L'(\l S_0,x,y\r,5,5)$, which is empty.
If $x$ is joined with $S_0$ by a double edge, we denote by $S_4\subset \l S_0,x\r$
a subdiagram of the type $B_4$ and consider the list $L'(\l S_0,x,y\r,5,5,S_4)$. The latter list consists of 
a unique diagram $\Sigma'$, shown in  Fig.\ref{5_h4_}(b). 
 
Let $S_4\subset  \Sigma'$ be the subdiagram of type $B_4$ marked by a gray box.
$S_4$ contains an end of the dotted edge and has a unique good neighbor (and no non-neighbors) in $\Sigma'$.
Hence, it has at least one good neighbor (or non-neighbor) in $\Sigma(P)\setminus \Sigma'$, 
so $\Sigma(P)$ contains a diagram from the list $L'(\Sigma',5,5,S_4)$, which is empty.

\begin{figure}[!h]
\begin{center}
\psfrag{a}{(a)}
\psfrag{b}{(b)}
\psfrag{6}{\small $6$}
\psfrag{6}{\small $8$}
\psfrag{x}{\small $x$}
\psfrag{y}{\small $y$}
\psfrag{s}{\small $\l S_0,x,y \r$}
\psfrag{s'}{\small $\Sigma'$}
\psfrag{3,4}{\tiny $3,\!4$}
\psfrag{3,4,5}{\tiny $3,\!4,\!5$}
\epsfig{file=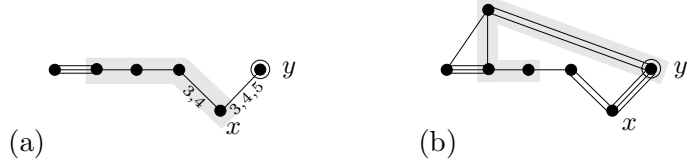,width=0.7\linewidth}
\caption{Notation to the proof of Lemma~\ref{5_h4}. (a) six possibilities for $\l S_0,x,y\r$; (b)  diagram $\Sigma'$.}
\label{5_h4_}
\end{center}
\end{figure}

\end{proof}

\begin{lemma}
\label{5_h3}
$\Sigma(P)$  contains no subdiagrams of the type $H_3$.

\end{lemma}

\begin{proof}
Suppose that $S_0\subset \Sigma(P)$ is a subdiagram of the type $H_3$. 
In view of Lemma~\ref{5_h4}, the diagram $S_0$ has no good neighbors, 
and $\o S_0 = \Sigma_{S_0}$ is a Lann\'er diagram of order 3 (see Lemma~\ref{polygon}).
By Lemmas~\ref{l'} and~\ref{5_mult}, $\Sigma(P)$ contains a subdiagram from the list $L'(\l S_0,\o S_0\r,5,5)$.
This list consists of 12 diagrams, 5 of which contain a subdiagram of the type $H_4$. Again, 
by Lemma~\ref{l'}, $\Sigma(P)$ contains a subdiagram from the list $L'(\Sigma_1,5,5)$, where $\Sigma_1$
ranges over the 7 diagrams of  $L'(\l S_0,\o S_0\r,5,5)$ containing no subdiagram of the type $H_4$.
All these lists   $L'(\Sigma_1,5,5)$ are empty, which completes the proof.

\end{proof}

\begin{lemma}
\label{5_h2}
$\Sigma(P)$  contains no subdiagrams of the type $G_2^{(5)}$.

\end{lemma}

\begin{figure}[!h]
\begin{center}
\psfrag{6}{\small $6$}
\psfrag{6}{\small $8$}
\psfrag{6}{\small $10$}
\psfrag{6}{\small $12$}
\epsfig{file=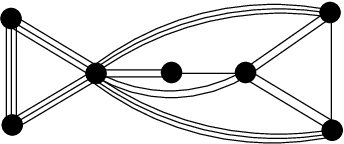,width=0.23\linewidth}
\caption{To the proof of Lemma~\ref{5_h2}.}
\label{5_h2_}
\end{center}
\end{figure}

\begin{proof}
Suppose that $S_0\subset \Sigma(P)$ is a subdiagram of the type $G_2^{(5)}$. 
Then $S_0$ has no good neighbors, and $\o S_0 = \Sigma_{S_0}$.
$P(S_0)$ is a 3-polytope with at most one pair of non-intersecting facets, 
so $\o S_0$ is either  is a Lann\'er diagram of order 4, 
or a diagram of a triangular prism. If $\o S_0$ is a diagram of a triangular prism, let  $\Sigma_1$ be a 
diagram spanned by $S_0$ and $\o S_0$ without tail.
In case of a Lann\'er diagram of order 4, let $\Sigma_1=\l\o S_0,S_0\r$.
By Lemmas~\ref{l'} and ~\ref{5_h3}, $\Sigma(P)$ contains a subdiagram from one of the lists $L'(\Sigma_1,5,5)$
with $\Sigma_1$ as above. Notice that we may consider only Lann\'er diagrams and diagrams of prisms 
not containing subdiagrams of the type $H_3$.
The union of these lists contains 5 entries, only one of them contains no subdiagram of the type $H_3$.
We present this diagram in  Fig.~\ref{5_h2_} and
denote it by $\Sigma_2$. By  Lemma~\ref{l'},  $\Sigma(P)$ contains a subdiagram from the list $L'(\Sigma_2,5,5)$,
which is empty.

\end{proof}

\begin{lemma}
\label{5_f4}
$\Sigma(P)$  contains no subdiagrams of the types $F_4$.

\end{lemma}

\begin{proof}
Suppose that $S_0\subset \Sigma(P)$ is a subdiagram of the type $F_4$.
Then $S_0$ has no good neighbors, so $\o S_0=\Sigma_{S_0}$ is a dotted edge.
Let $S_1\subset S_0$ be a subdiagram of the type $B_3$.
$P(S_1)$ is a 2-polytope with a pair of non-intersecting facets, 
so $\Sigma(P)$ contains a node $x$ such that $x$ is not a bad neighbor of $S_1$, and
the edge $xt_1$ turns into a dotted edge in $\Sigma_{S_1}$.
It follows from~\cite[Theorem~2.2]{All}  that $\l S_0,x\r$ is one of the two diagrams $\Sigma_1$ and $\Sigma_2$ 
shown in   Fig.~\ref{5_f4_}(a).
Notice, that $x$ is a bad neighbor of $S_0$, so it is joined with at least one end (denote it by $y$) of the dotted edge 
(otherwise the diagram $\l S_0,x,\o S_0 \r$ is superhyperbolic). By Lemmas~\ref{5_mult} and~\ref{5_h2}, $[y,x]=3$ or 4.
In case of the diagram $\Sigma_1$ this leads to a parabolic subdiagram of the type $\w F_4$ or $\w C_3$.
In case of  $\Sigma_2$ this implies that  $[y,x]=3$ (otherwise we obtain a parabolic subdiagram of the type $\w C_4$).
So, we are left with the only possibility for the diagram $\l \Sigma_2,x\r$, see Fig.~\ref{5_f4_}(b).
By Lemma~\ref{l'}, $\Sigma(P)$ contains a subdiagram from the list $L'(\l \Sigma_2,x\r,4,5)$.
However, this list is empty.

\begin{figure}[!h]
\begin{center}
\psfrag{a}{(a)}
\psfrag{b}{(b)}
\psfrag{x}{\small $x$}
\psfrag{y}{\small $y$}
\psfrag{s1}{\small $\Sigma_1$}
\psfrag{s2}{\small $\Sigma_2$}
\epsfig{file=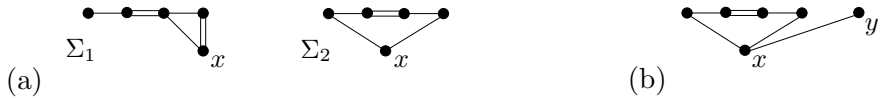,width=0.9\linewidth}
\caption{To the proof of Lemma~\ref{5_f4}.}
\label{5_f4_}
\end{center}
\end{figure}

\end{proof}

\begin{lemma}
\label{5_b5}
$\Sigma(P)$  contains no subdiagrams of the type $B_5$.

\end{lemma}

\begin{proof}
Suppose that $S_0\subset \Sigma(P)$ is a subdiagram of the type $B_5$. 
Let $S_1\subset S_0$ be a subdiagram of the type $B_4$. $P(S_0)$ is a 1-polytope,
so $\Sigma_{S_0}$ is a dotted edge. By~\cite[Theorem~2.2]{All}, this may happen only if $\l S_1,\o S_1 \r$ is
one of two diagrams $\Sigma_1$ and $\Sigma_2$ shown in the left row of
Table~\ref{5_b5_}. 

\begin{table}[!h]
\begin{center}
\caption{Notation to the proof of Lemma~\ref{5_b5}.}
\label{5_b5_}
\vspace{10pt}
\psfrag{1}{\scriptsize $t_1$}
\psfrag{2}{\scriptsize $t_2$}
\psfrag{3}{\scriptsize $t_3$}
\psfrag{4}{\scriptsize $t_4$}
\psfrag{5}{\scriptsize $t_5$}
\psfrag{L'1}{\tiny $L'(\Sigma_1,\!4,\!5)$}
\psfrag{L'2}{\tiny $L'(\Sigma_1^{1a},\!4,\!5)$}
\psfrag{L'3}{\tiny $L'(\Sigma_2^1,\!4,\!5)$}
\psfrag{s1}{\small $\Sigma_1$}
\psfrag{s2}{\small $\Sigma_2$}
\psfrag{s1a}{\small $\Sigma_1^{1a}$}
\psfrag{s1b}{\small $\Sigma_1^{1b}$}
\psfrag{s21}{\small $\Sigma_2^1$}
\psfrag{s22}{\small $\Sigma_2^2$}
\psfrag{s12}{\small $\Sigma_1^2$}
\psfrag{x}{\scriptsize $x$}
\psfrag{2,3}{\tiny $2,\!3$}
\epsfig{file=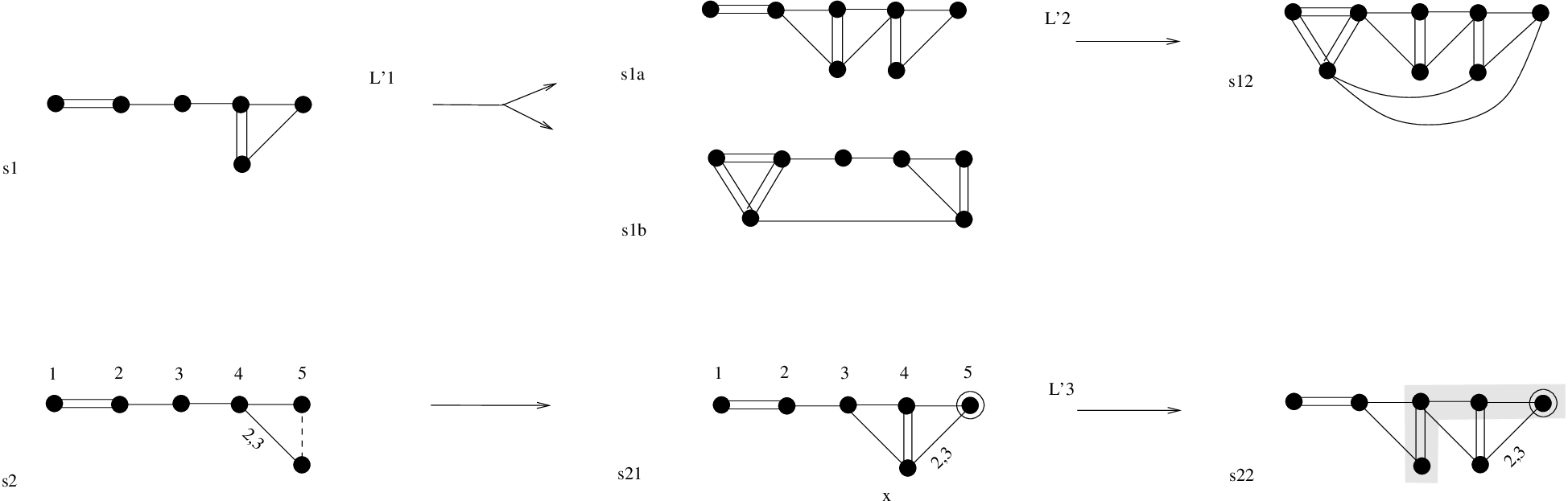,width=0.997\linewidth}
\end{center}
\end{table}

Consider the diagram  $\Sigma_1$. By Lemma~\ref{l'}, $\Sigma(P)$ contains a diagram from the list $L'(\Sigma_1,4,5)$ 
The list consists of two diagrams $\Sigma_1^{1a}$ and  $\Sigma_1^{1b}$
(see Table~\ref{5_b5_}). The diagram $\Sigma_1^{1b}$ contains a subdiagram of the type $F_4$,
which is impossible by Lemma~\ref{5_f4}. For the diagram $\Sigma_1^{1a}$ we consider the list $L'(\Sigma_1^{1a},4,5)$,
which consists of a unique diagram  $\Sigma_1^{2}$. The latter diagram contains a subdiagram of the type $F_4$,
which is impossible.

Consider the diagram  $\Sigma_2$.  Let $S_2\subset \Sigma_2$ be a subdiagram of the type $B_3$.
$P(S_2)$ is a polygon with at least 4 edges. So, there exists at least one good neighbor or a non-neighbor $x$ of $S_2$
such that $xt_4$ turns into a dotted edge in $\Sigma_{S_2}$ (see Table~\ref{5_b5_} for the notation).
This is possible only if $[x,t_3]=3$ and   $[x,t_4]=4$.
Notice, that  $[x,t_5]\ne 4$, otherwise $\l S_0,x\r$ contains a parabolic subdiagram of the type $\w C_4$.
Denote by $\Sigma_2^1$ the subdiagram $\l S_0,x\r$  (see Table~\ref{5_b5_}).
By Lemma~\ref{l'}, $\Sigma(P)$ contains a diagram from the list $L'(\Sigma_2^1,4,5)$,
which consists of a unique diagram $\Sigma_2^2$   (see Table~\ref{5_b5_} again).
Consider the subdiagram $S_3\subset\Sigma_2^2$ marked by a gray box.
$S_3$ is a diagram of the type $B_4$ containing an end $t_5$ of the dotted edge.
So, $\l S_3,\o S_3 \r$ is a diagram of the same type as $\Sigma_1$.
As it is shown above, the diagram  $\Sigma_1$ cannot be a subdiagram of $\Sigma(P)$.
So, the diagram $S_3$ also cannot be a subdiagram of $\Sigma(P)$, which completes the proof.

\end{proof}

\begin{lemma}
\label{d5}
The Main Theorem holds in dimension 5. 

\end{lemma}

\begin{proof}
Let $P$ be a compact hyperbolic Coxeter 5-polytope with at least 5 facets and exactly one pair of non-intersecting facets.
By Lemmas~\ref{5_mult}-\ref{5_b5}, $\Sigma(P)$ does not contain neither edges of multiplicity greater than 2, 
nor diagrams of the type $B_5$. Applying Lemmas~\ref{b_k} and~\ref{b_2}, we finish the proof.   

\end{proof}

\medskip
\noindent
{\bf Remark.} Instead of Lemmas~\ref{5_h4}-\ref{5_b5} one could use the reasoning similar to the proof of Lemma~\ref{4,b};
however, in dimension 5 this leads to very long computation (in particular, one should find the list 
$L'(\l T_1,T_2 \r,5,5)$,
where $T_1$ and $T_2$ are Lann\'er diagrams of order 3 containing no multi-multiple edges,
and then for each diagram $\Sigma\in  L'(\l T_1,T_2 \r,5,5)$ we should find the list  $L'(\Sigma,5,5)$).

\subsection{Dimension 6}

Let $P$ be a 6-dimensional compact hyperbolic Coxeter polytope such that $\Sigma(P)$ contains a unique dotted edge and 
$P$ has at least 10 facets.

\begin{lemma}
\label{6_mult}
$\Sigma(P)$  contains no multi-multiple edges.

\end{lemma}

\begin{proof}
Suppose that $S_0\subset \Sigma(P)$ is a multi-multiple edge of the maximum multiplicity in $\Sigma(P)$. 
Then $S_0$ has no good neighbors, and, by Lemma~\ref{l1}, $\Sigma(P)$ contains a subdiagram $\l S_0,y_1,y_0,S_1 \r$ from 
the list $L_1(6)$. The list consists of 8  diagrams shown in Table~\ref{6_multi}. 
We denote these diagrams $\Sigma_1,\dots,\Sigma_8$.
Notice, that for each of the diagrams it is easy to find out where the subdiagram $S_0$ is 
(the multi-multiple edge with a unique bad neighbor), where the node $y_1$ is (which is the bad neighbor of $S_0$), and
where $\l y_0,S_1\r$ is. The node $y_1$ is a bad neighbor 
of the subdiagram $S\subset\l y_0,S_1\r $ of the type $H_4$
or $F_4$, so the node $\l y_0,S_1\r\setminus S$ is an end of the dotted edge (we mark the end of the dotted edge by
a circle). For each of $\Sigma_1,\dots,\Sigma_8$ (except $\Sigma_7$) denote by $S_2$ the 
elliptic subdiagram of order 5 marked by 
a gray box. Notice, that $S_2$ has a unique good neighbor (or a unique non-neighbor) in $\Sigma_i$.
So, it has one more in $\Sigma(P)$. Thus, in case of diagrams $\Sigma_1,\dots,\Sigma_6$ we consider the lists
$L'(\Sigma_i,k(\Sigma_i),6,S_2)$, where $k(\Sigma_i)=6$ for $i=1,2,3$ and  $k(\Sigma_i)=10$ for $i=4,5,6$.
The lists are empty. 

\begin{table}[!h]
\begin{center}
\caption{The list $L_1(6)$.}
\label{6_multi}
\vspace{10pt}
\psfrag{a}{(a)}
\psfrag{b}{(b)}
\psfrag{6}{\scriptsize $6$}
\psfrag{8}{\scriptsize $8$}
\psfrag{10}{\scriptsize $10$}
\psfrag{s1}{\small $\Sigma_1$}
\psfrag{s2}{\small $\Sigma_2$}
\psfrag{s3}{\small $\Sigma_3$}
\psfrag{s4}{\small $\Sigma_4$}
\psfrag{s5}{\small $\Sigma_5$}
\psfrag{s6}{\small $\Sigma_6$}
\psfrag{s7}{\small $\Sigma_7$}
\psfrag{s8}{\small $\Sigma_8$}
\epsfig{file=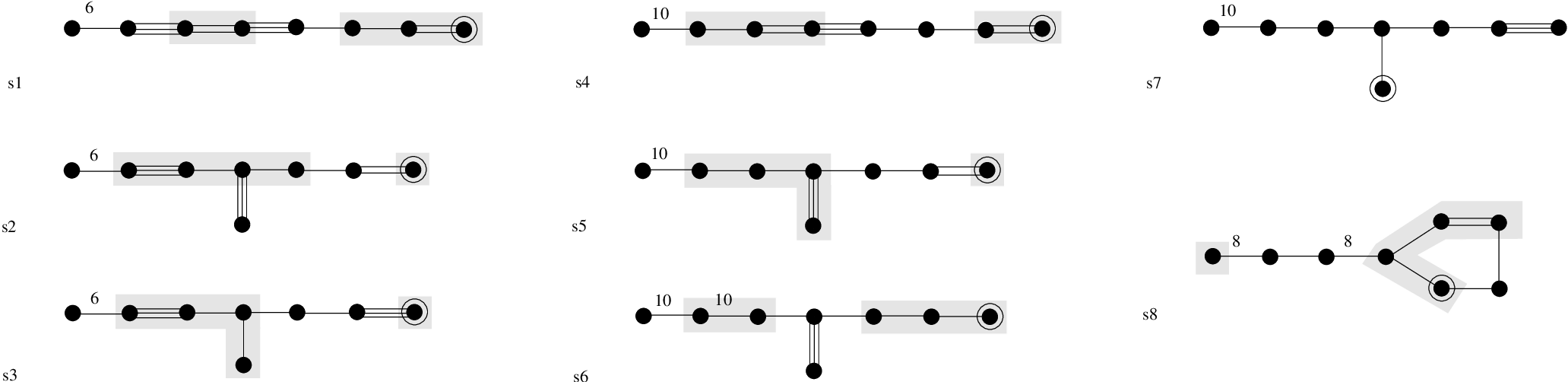,width=0.997\linewidth}
\end{center}
\end{table}

We are left to consider the diagrams $\Sigma_7$ and $\Sigma_8$.
In case of the diagram $\Sigma_7$ denote by $\Sigma_7^1\subset \Sigma_7$ the subdiagram with 
the end of the dotted edge discarded.
Let $S_2\subset \Sigma_7^1$ be the subdiagram of the type $H_4$. Since $S_2$ has only two non-neighbors in $\Sigma_7^1$, 
it has at least one more in $\Sigma(P)$. So,   $\Sigma(P)$ contains a diagram from the list $L'(\Sigma_7^1,10,6,S_2^{})$,
which consists of two diagrams shown in Fig.~\ref{6_multi2}.
The diagram shown in Fig.~\ref{6_multi2}(a) is a  diagram of a 6-polytope with 9 facets,
so by~\cite[Lemma~1]{nodots} it cannot be a subdiagram of $\Sigma(P)$.  
Denote by $\Sigma_7^{2}$ the diagram shown in Fig.~\ref{6_multi2}(b) and consider
the elliptic subdiagram $S_3\subset \Sigma_7^{2}$ of order 5 marked by a gray box. 
It has no good neighbors (non-neighbors) in  $\Sigma_7^{2}$,
so at least one of its good neighbors (non-neighbors) is not joined with    $\Sigma_7^{2}$ by a dotted edge.
However, the list   $L'(\Sigma_7^{2},10,6,S_3^{})$ is empty, and 
the diagram  $\Sigma_7$ cannot be a subdiagram of $\Sigma(P)$.
 
\begin{figure}[!h]
\begin{center}
\psfrag{a}{(a)}
\psfrag{b}{(b)}
\psfrag{6}{\scriptsize $6$}
\psfrag{8}{\scriptsize $8$}
\psfrag{10}{\scriptsize $10$}
\epsfig{file=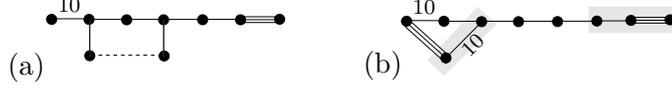,width=0.7\linewidth}
\caption{Treating the diagram $\Sigma_7$, see Lemma~\ref{6_mult}.}
\label{6_multi2}
\end{center}
\end{figure}

Consider the remaining diagram, $\Sigma_8$.
The subdiagram $S_2$ of order 5 (marked by a gray box) 
has a unique good neighbor in $\Sigma_8$. $S_2$ contains an end of the dotted edge, so,
the second good neighbor of $S_2$ (or non-neighbor) is not joined with $\Sigma_8$ by the dotted edge. 
Therefore, $\Sigma(P)$ contains a diagram from the 
list  $L'(\Sigma_8,8,6,S_2)$, which consists of a unique diagram $\Sigma_8^1$ shown in Table~\ref{6_multi3}.
Let $S_3\subset \Sigma_8^1$ be a subdiagram of order 4 marked by a gray box (see Table~\ref{6_multi3}). $S_3$ 
has only one non-neighbor (and no good neighbors) in $\Sigma_8^1$,
so it should have at least two more in $\Sigma(P)$. Therefore,  $\Sigma(P)$ contains a diagram from the list 
 $L'(\Sigma_8^{1},8,6,S_3)$, which consists of two diagrams $\Sigma_8^{2a}$ and $\Sigma_8^{2b}$ 
shown in Table~\ref{6_multi3}.
Denote by $\Sigma_8^{2a'}$ and $\Sigma_8^{2b'}$ these diagrams with the end of the dotted edge discarded.
Denote by $S_4$ the subdiagram of order 4 in  $\Sigma_8^{2a'}$ and $\Sigma_8^{2b'}$ marked by a gray box. 
$S_4$ has only to non-neighbors (and no good neighbors) in $\Sigma_8^{2a'}$  (and in $\Sigma_8^{2b'}$), 
so, it has at least one more in $\Sigma(P)$.
Since the diagrams $\Sigma_8^{2a'}$ and $\Sigma_8^{2b'}$ contain no end of dotted edge,
$\Sigma(P)$  contains a diagram from one of the lists $ L'(\Sigma_8^{2a'},8,6,S_4)$ 
and $ L'(\Sigma_8^{2b'},8,6,S_4)$.
The first of these lists is empty, the second one consists of two diagrams $\Sigma_8^{3a'}$ and $\Sigma_8^{3b'}$ 
shown  in Table~\ref{6_multi3}. Returning the end of the dotted edge and computing the weight of the edge joining
that with $\Sigma_8^{3a'}\setminus\Sigma_8^{2a'}$ (resp., with $\Sigma_8^{3b'}\setminus\Sigma_8^{2b'}$), 
we obtain subdiagrams $\Sigma_8^{3a}$ and $\Sigma_8^{3b}$ of $\Sigma(P)$, see Table~\ref{6_multi3}.

Consider the diagram  $\Sigma_8^{3a}$. 
Let $S_5\subset \Sigma_8^{3a}$ be a subdiagram of the type $D_4$ marked by a gray box. It has only two non-neighbors 
(and no good neighbors) in $\Sigma_8^{3a}$. Hence, $\Sigma_8^{3a}$ is not a diagram of a Coxeter polytope. 
Now, consider the diagram  $\Sigma_8^{3a'}$.
Since there exists a good neighbor (or a non-neighbor) of $S_5$ which does not belong to  $\Sigma_8^{3a}$,
we conclude that $\Sigma(P)$ contains a diagram from the list  $ L'(\Sigma_8^{3a'},8,6,S_5)$, which is empty.

We are left to consider the diagram  $\Sigma_8^{3b}$. 
Consider the diagram $S_6$ of the type $G_2^{(8)}$ marked by a gray box. 
It has no good neighbors in $\Sigma(P)$, so $\o S_6=\Sigma_{S_6}$
is either a Lann\'er diagram of order 5 or an Esselmann diagram (since one of the ends of the dotted edge is
a bad neighbor of $S_6$).
However, discarding from     $\Sigma_8^{3b}$ the subdiagram $S_6$ with all its bad neighbors, 
we obtain a subdiagram $\Sigma'$ shown in 
Table~\ref{6_multi3}, which is neither a Lann\'er diagram nor a part of an Esselmann diagram.
Therefore,  the diagram  $\Sigma_8$ also cannot be a subdiagram of $\Sigma(P)$, and the lemma is proved.

\begin{table}[!h]
\begin{center}
\caption{Treating the diagram $\Sigma_8$, see Lemma~\ref{6_mult}.}
\label{6_multi3}
\vspace{10pt}
\psfrag{a}{$\Sigma_8^1$}
\psfrag{b}{$\Sigma_8^{2a}$}
\psfrag{c}{$\Sigma_8^{2b}$}
\psfrag{d}{$\Sigma_8^{3a}$}
\psfrag{e}{$\Sigma_8^{3b}$}
\psfrag{f}{$\Sigma'$}
\psfrag{g}{$\Sigma_8^{3a'}$}
\psfrag{h}{$\Sigma_8^{3b'}$}
\psfrag{6}{\small $6$}
\psfrag{8}{\tiny $8$}
\psfrag{10}{\small $10$}
\epsfig{file=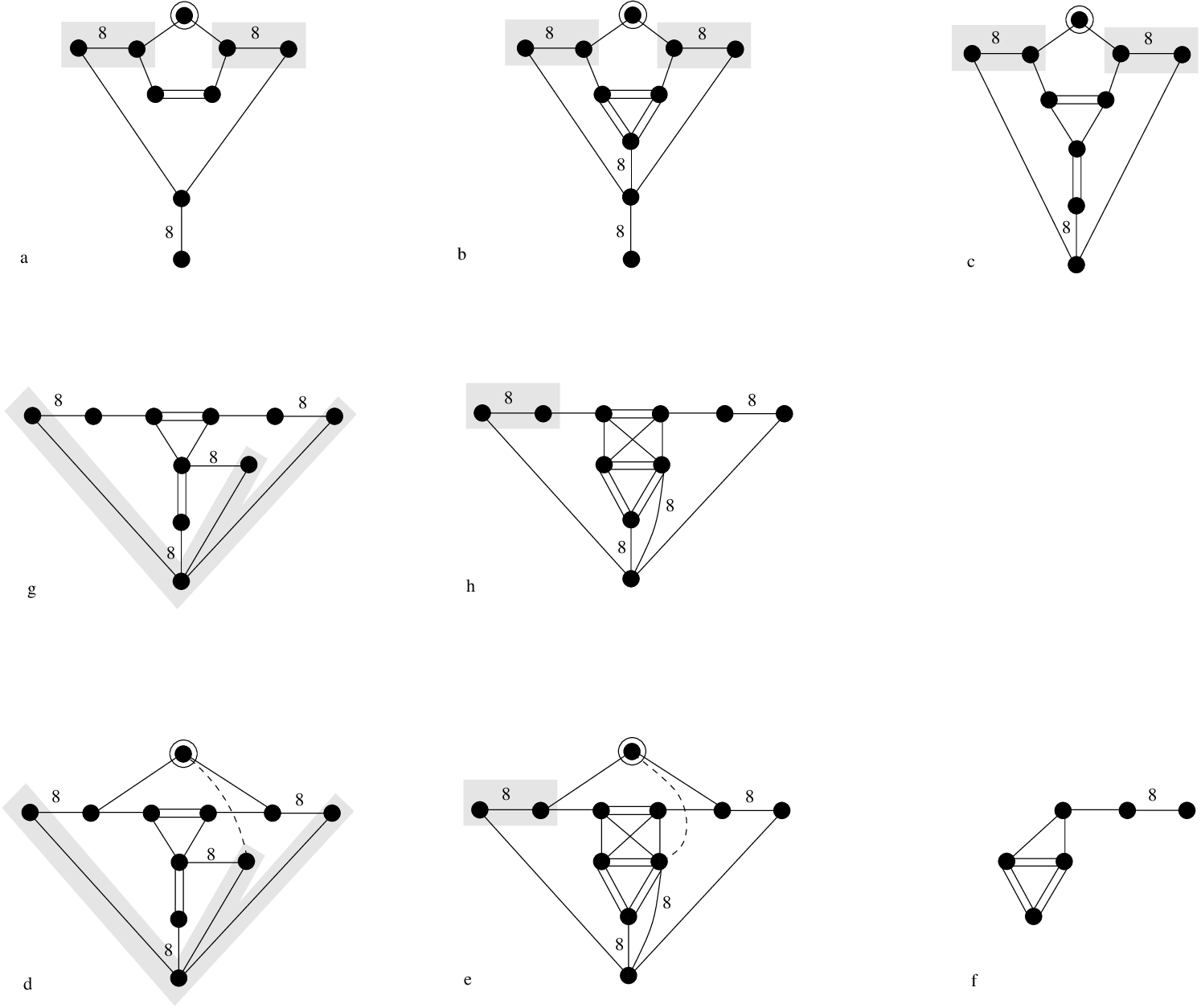,width=0.9997\linewidth}
\end{center}
\end{table}

\end{proof}

\begin{lemma}
\label{6_h4,f4}
$\Sigma(P)$  contains no subdiagrams of the types $H_4$ and $F_4$.

\end{lemma}

\begin{proof}
Suppose that $S_0\subset \Sigma(P)$ is a subdiagram of the type $H_4$ or $F_4$.
Then $\Sigma(P)$ contains a diagram from the list $L_1(H_4,6)$ or $L_1(F_4,6)$.
The union of these lists consists of 9 diagrams  shown in Table~\ref{6_h4_f4}, we denote these diagrams 
$\Sigma_1,\dots,\Sigma_9$ (the list  $L_1(H_4,6)$ is shown in the left column,  $L_1(F_4,6)$ is shown in the right one). 
For the diagrams $\Sigma_1,\dots,\Sigma_6$ we consider the lists 
$L'(\Sigma_i,5,6)$, which turn out to be empty. In particular, this implies that 
$\Sigma(P)$  contains no subdiagram of the type $H_4$.

\begin{table}[!h]
\begin{center}
\caption{Lists $L_1(H_4,6)$ and $L_1(F_4,6)$.}
\label{6_h4_f4}
\vspace{10pt}
\psfrag{a}{(a)}
\psfrag{b}{(b)}
\psfrag{6}{\small $6$}
\psfrag{8}{\small $8$}
\psfrag{10}{\small $10$}
\psfrag{s1}{\small $\Sigma_1$}
\psfrag{s2}{\small $\Sigma_2$}
\psfrag{s3}{\small $\Sigma_3$}
\psfrag{s4}{\small $\Sigma_4$}
\psfrag{s5}{\small $\Sigma_5$}
\psfrag{s6}{\small $\Sigma_6$}
\psfrag{s7}{\small $\Sigma_7$}
\psfrag{s8}{\small $\Sigma_8$}
\psfrag{s9}{\small $\Sigma_9$}
\epsfig{file=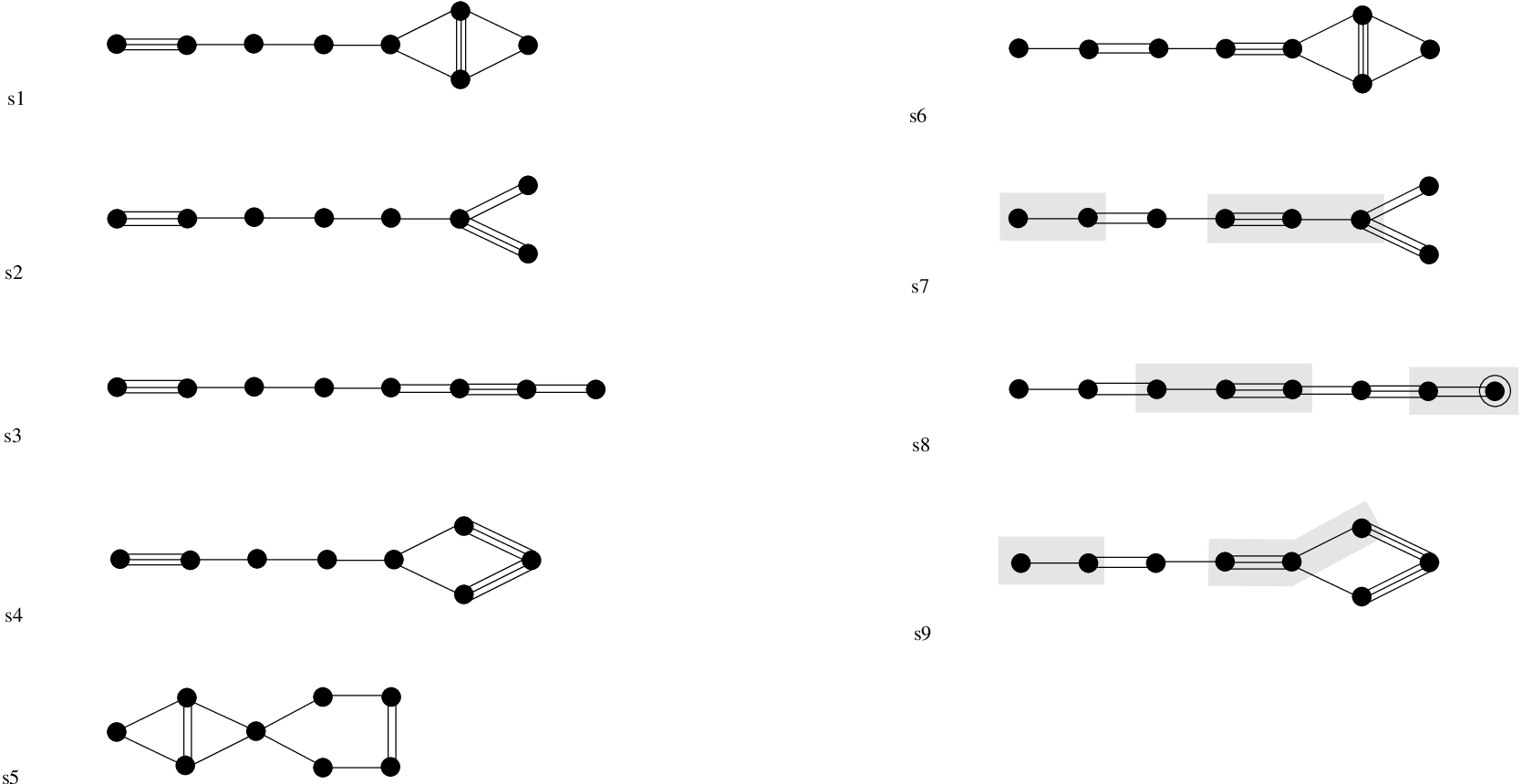,width=0.8\linewidth}
\end{center}
\end{table}

For the diagrams $\Sigma_7$, $\Sigma_8$ and $\Sigma_9$ we denote by $S_2$ a subdiagram of order $5$ marked by a gray box.
It has neither good neighbors nor non-neighbors in cases of $\Sigma_7 $ and $\Sigma_9$, and it has a unique 
good neighbor in case of  $\Sigma_8$, however in the latter case $S_2$ contains an end of the dotted edge 
(we know where the end of the dotted edge is, 
since $y_1$ is a good neighbor of a subdiagram of the type $B_2\subset \o S_0$ 
but not of the subdiagram of the type $G_2^{(5)}$, which is maximal). 
Therefore, $\Sigma(P)$ contains a subdiagram from one of the lists $L'(\Sigma_i,5,6,S_2)$, $i=7,8,9$.
Each of the lists $L'(\Sigma_7,5,6,S_2)$ and $L'(\Sigma_8,5,6,S_2)$ consist of the diagram $\Sigma^{78}$ 
shown in Fig.~\ref{6_h4f4_2}(a),
 the list $L'(\Sigma_9,5,6,S_2)$ consists of the diagram $\Sigma^{9}$  shown in Fig.~\ref{6_h4f4_2}(b).
For each of  $\Sigma^{78}$ and $\Sigma^{9}$ consider a subdiagram $S_3$ of the type $H_3$ marked by a gray box. 
As it was shown above, $S_3$ has no good neighbors in $\Sigma(P)$.
So, $P(S_3)$ is a 3-polytope with at most one pair of non-intersecting facets, and 
$\o S_3=\Sigma_{S_3}$ is either a Lann\'er diagram of order 4, or a diagram of a 3-prism. 
The former case is impossible since $\o S_3$ contains a Lann\'er subdiagram of order 3, so $P(S_3)$
is a prism.
In case of the diagram $\Sigma^9$ this implies that $S_3$ has at least 2 additional non-neighbors, and hence,
$\Sigma(P)$ contains a diagram from the list $L'(\Sigma^9,5,6,S_3)$, which is empty.

We are left with the diagram  $\Sigma^{78}$. 
Let $T$ be the Lann\'er subdiagram of $\Sigma^{78}$ contained in $\o S_3$, and let $x$ be  
the leaf of   $\Sigma^{78}$ (node of valency 1). 
Since $P(S_3)$ is a prism,
there exists a non-neighbor of $S_3$, a node $y\in \Sigma(P)\setminus \Sigma^{78}$, 
such that $y$ is joined with $T$ by some edge and $y$ is joined with $x$ by a dotted edge.
However, the list $L'(\Sigma^{78}\setminus x,5,6,S_3)$ contains no entry in which the new node is joined with $T$.
This completes the proof.

\begin{figure}[!h]
\begin{center}
\psfrag{a}{(a)}
\psfrag{b}{(b)}
\psfrag{6}{\small $6$}
\psfrag{8}{\small $8$}
\psfrag{10}{\small $10$}
\psfrag{s78}{\small $\Sigma^{78}$}
\psfrag{s9}{\small $\Sigma^9$}
\epsfig{file=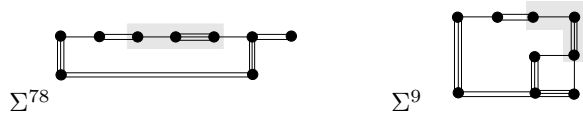,width=0.6\linewidth}
\caption{Treating the diagrams $\Sigma_7$, $\Sigma_8$ and $\Sigma_9$, see Lemma~\ref{6_h4,f4}.}
\label{6_h4f4_2}
\end{center}
\end{figure}

\end{proof}

\begin{lemma}
\label{6_h3}
$\Sigma(P)$  contains no subdiagram of the type $H_3$.

\end{lemma}

\begin{proof}
Suppose that $S_0\subset \Sigma(P)$ is a subdiagram of the type $H_3$.
Then $\Sigma(P)$ contains a diagram from the list $L_1(H_3,6)$, which consists of 4 diagrams. 
Two of these diagrams contain the subdiagrams of the type 
$F_4$ or $H_4$. The remaining two diagrams are the diagrams $\Sigma_1$ and $\Sigma_2$ shown in  Fig.~\ref{6_h3_}.
For the diagram   $\Sigma_1$ we check the list $L'(\Sigma_1,5,6)$, which is empty.
For the diagram  $\Sigma_2$ the list $L'(\Sigma_2,5,6)$ consists of a unique entry  $\Sigma_2'$ (see Fig.~\ref{6_h3_}).
Let $S_2\subset \Sigma_2'$ be a subdiagram of the type $B_2$ marked by a gray box.
Discarding from  $ \Sigma_2'$ the subdiagram $S_2$ with all its bad neighbor, we obtain
a subdiagram $\Theta$ of order 5 which consists of a Lann\'er diagram of order 3 and of two separate nodes. 
It is easy to see that $\Theta $ is not a subdiagram of a Lann\'er diagram of order 5, 
of an Esselmann diagram or of diagram of a 4-prism.
Therefore, $\Sigma_{S_2}$ contains at least 7 nodes, and $\Sigma(P)$ 
contains  a diagram from the list $L'(\Sigma_2',5,6,S_2)$, which is empty.

\begin{figure}[!h]
\begin{center}
\psfrag{a}{(a)}
\psfrag{b}{(b)}
\psfrag{6}{\small $6$}
\psfrag{8}{\small $8$}
\psfrag{10}{\small $10$}
\psfrag{s1}{\small $\Sigma_{1}$}
\psfrag{s2}{\small $\Sigma_2$}
\psfrag{s3}{\small $\Sigma_2'$}
\epsfig{file=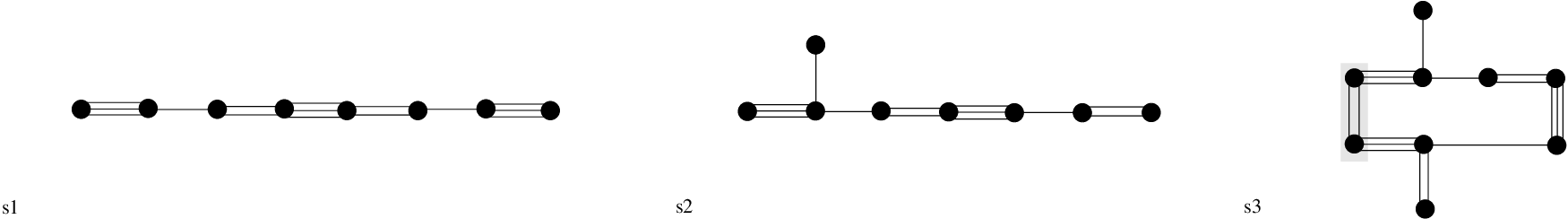,width=0.95\linewidth}
\caption{To the proof of Lemma~\ref{6_h3}, see Lemma~\ref{6_h3}.}
\label{6_h3_}
\end{center}
\end{figure}

\end{proof}

\begin{lemma}
\label{6_h2}
$\Sigma(P)$  contains no subdiagram of the type $G_2^{(5)}$.

\end{lemma}

\begin{proof}
Suppose that $S_0\subset \Sigma(P)$ is a subdiagram of the type $G_2^{(5)}$.
Then $S_0$ has no good neighbors, so $P(S_0)$ is a 4-polytope with at most one pair of non-intersecting facets, 
so (by the Main Theorem in dimension $d=4$),
a 4-polytope with at most $7$ facets. There are only four 4-polytopes with at most $7$ facets such that 
their Coxeter diagrams contain no subdiagram of the type $H_4$ or $F_4$. The diagrams are shown in Fig.~\ref{6_h2_}(a)
(the diagram $\Sigma_1$ corresponds to two 4-prisms).
Notice, that all these diagrams contain dotted edges. At the same time, the diagram $\Sigma_3$ contains a subdiagram 
$S_1$ of the type 
$G_2^{(5)}$ such that $\o S_1$ definitely contains no dotted edges (one end of the dotted edge is a bad neighbor of $S_1$).
This is impossible, so we are left with the diagrams $\Sigma_1$ and $\Sigma_2$. 
Denote by $\Sigma_1'$ and $\Sigma_2'$ the diagrams with respectively one and two nodes discarded (see  Fig.~\ref{6_h2_}(b)).
Let $S_2$ be a subdiagram of  $\Sigma_1'$ or $\Sigma_2'$ of the type $B_4$ (marked by a gray box).
The diagram $S_2$ has only two good neighbors in $\l S_0,\Sigma_1'\r$ as well as in $\l S_0,\Sigma_2'\r$,
at the same time, $S_2$ contains an end of the dotted edge. Therefore, $S_2$ has a good neighbor (or a non-neighbor) in
$\Sigma(P)\setminus  \l S_0,\Sigma_1'\r$ (or in $\Sigma(P)\setminus  \l S_0,\Sigma_2'\r$ respectively),
and $\Sigma(P)$ contains a diagram from the list $L'(\l S_0,\Sigma_1'\r,5,6,S_2)$ or  $L'(\l S_0,\Sigma_2'\r,5,6,S_2)$.
Both these lists are empty, and the lemma is proved.

\begin{figure}[!h]
\begin{center}
\psfrag{a}{(a)}
\psfrag{b}{(b)}
\psfrag{6}{\small $6$}
\psfrag{8}{\small $8$}
\psfrag{2,3}{\tiny $2,\!3$}
\psfrag{s1}{\small $\Sigma_1$}
\psfrag{s2}{\small $\Sigma_2$}
\psfrag{s3}{\small $\Sigma_3$}
\psfrag{s4}{\small $\Sigma_1'$}
\psfrag{s5}{\small $\Sigma_2'$}
\epsfig{file=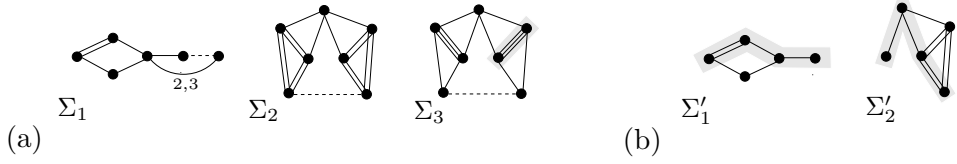,width=0.99\linewidth}
\caption{To the proof of Lemma~\ref{6_h2}. 
(a) 4-polytopes with at most $7$ facets containing no subdiagrams $H_4$, $F_4$ and $G_2^{(k)}$, $k\ge 6$; 
(b) some subdiagrams of the diagrams shown in (a) ($\Sigma_1'\subset \Sigma_1$, $\Sigma_2'\subset \Sigma_2$).}
\label{6_h2_}
\end{center}
\end{figure}

\end{proof}

\begin{lemma}
\label{d6}
The Main Theorem holds in dimension 6. 

\end{lemma}

\begin{proof}
Let $P$ be a compact hyperbolic Coxeter 6-polytope with at least 10 facets and exactly one pair 
of non-intersecting facets.
By Lemmas~\ref{6_mult}-\ref{6_h2}, $\Sigma(P)$ does not contain edges of multiplicity greater than 2. 
Now we apply Lemmas~\ref{b_d}, and~\ref{b_2} to complete the proof.   

\end{proof}

\subsection{Dimension 7}

Let $P$ be a 7-dimensional hyperbolic Coxeter polytope such that $\Sigma(P)$ contains a unique dotted edge and 
$P$ has at least 11 facets.

\begin{lemma}
\label{7_mult}
$\Sigma(P)$  contains no multi-multiple edges.

\end{lemma}

\begin{proof}
Suppose that $S_0\subset \Sigma(P)$ is a multi-multiple edge of the maximum multiplicity in $\Sigma(P)$. 
Then $S_0$ has no good neighbors and $P(S_0)$ is either a 5-prism or a 5-polytope with 8 facets with a unique pair 
of non-intersecting facets (there is a unique such polytope). 
By Lemma~\ref{l1}, $\Sigma(P)$ contains a subdiagram $\l S_0,y_1,y_0,S_1 \r$ from 
the list $L_1(7)$. The list consists of 5  diagrams $\Sigma_1,\dots,\Sigma_5$ (see Table~\ref{7_multi}).
Notice, that for each of these diagrams the subdiagram $\l y_0,S_1\r$ is a part of a diagram of a 5-prism,
and we know where the end of the dotted edge is.
Denote by $S_2\subset \Sigma_i$, $i=1,\dots,5$ the elliptic subdiagram of order 6 marked by a gray box.
The diagram $S_2$ contains an end of the dotted edge and has at most 1 good neighbor in $\Sigma_i$.
Therefore, there exists a good neighbor or a non-neighbor of $S_2$ which is not joined with $\Sigma_i$ 
by a dotted edge. So, $\Sigma(P)$ contains a subdiagram from the list $L'(\Sigma_i,k(\Sigma_i),7)$,
where $\Sigma_i$ ranges over 5 diagrams  $\Sigma_1,\dots,\Sigma_5$ and $k(\Sigma_i)$ is a maximum multiplicity of
the edge in $\Sigma_i$. All these lists are empty, and the lemma is proved.

\begin{table}[!h]
\begin{center}
\caption{The list $L_1(7)$.}
\label{7_multi}
\vspace{10pt}
\psfrag{a}{(a)}
\psfrag{b}{(b)}
\psfrag{6}{\scriptsize $6$}
\psfrag{8}{\scriptsize $8$}
\psfrag{10}{\scriptsize $10$}
\epsfig{file=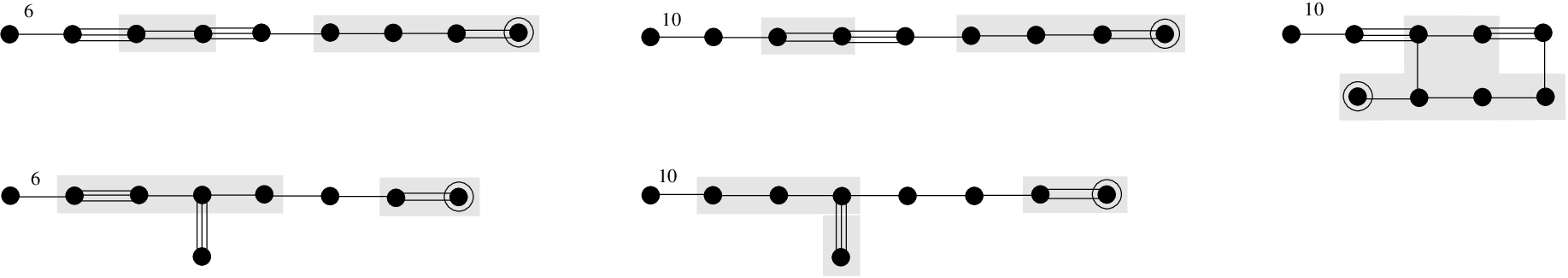,width=0.99\linewidth}
\end{center}
\end{table}

\end{proof}

\begin{lemma}
\label{7_h4,f4}
$\Sigma(P)$  contains no subdiagrams of the types $H_4$ and $F_4$.

\end{lemma}

\begin{proof}
Suppose that $S_0\subset \Sigma(P)$ is a subdiagram of the type $H_4$ or $F_4$.
Then $\Sigma(P)$ contains a diagram from the list $L_1(H_4,7)$ or $L_1(F_4,7)$.
Each of these lists consists of 3 diagrams, we denote these 6 diagrams by $\Sigma_1,\dots,\Sigma_6$ 
(see  Table~\ref{7_f4_h4}).
Notice that in cases of the diagrams   $\Sigma_2$, $\Sigma_3$, $\Sigma_5$ and $\Sigma_6$
we know where the end of the dotted edge is, since $y_1$ (the bad neighbor of $S_0$) is a good neighbor of a 
diagram $S_1\subset \o S_0$
of the type $B_3$, but not $H_3$.

\begin{table}[!hb]
\begin{center}
\caption{The lists $L_1(H_4,7)$ and  $L_1(F_4,7)$.}
\label{7_f4_h4}
\vspace{10pt}
\psfrag{s1}{$\Sigma_1$}
\psfrag{s2}{$\Sigma_2$}
\psfrag{s3}{$\Sigma_3$}
\psfrag{s4}{$\Sigma_4$}
\psfrag{s5}{$\Sigma_5$}
\psfrag{s6}{$\Sigma_6$}
\psfrag{s1}{$\Sigma_1$}
\psfrag{s2}{$\Sigma_2$}
\psfrag{s3}{$\Sigma_3$}
\psfrag{s4}{$\Sigma_4$}
\psfrag{s5}{$\Sigma_5$}
\psfrag{s6}{$\Sigma_6$}
\psfrag{1}{\scriptsize $t_1$}
\psfrag{2}{\scriptsize $t_2$}
\epsfig{file=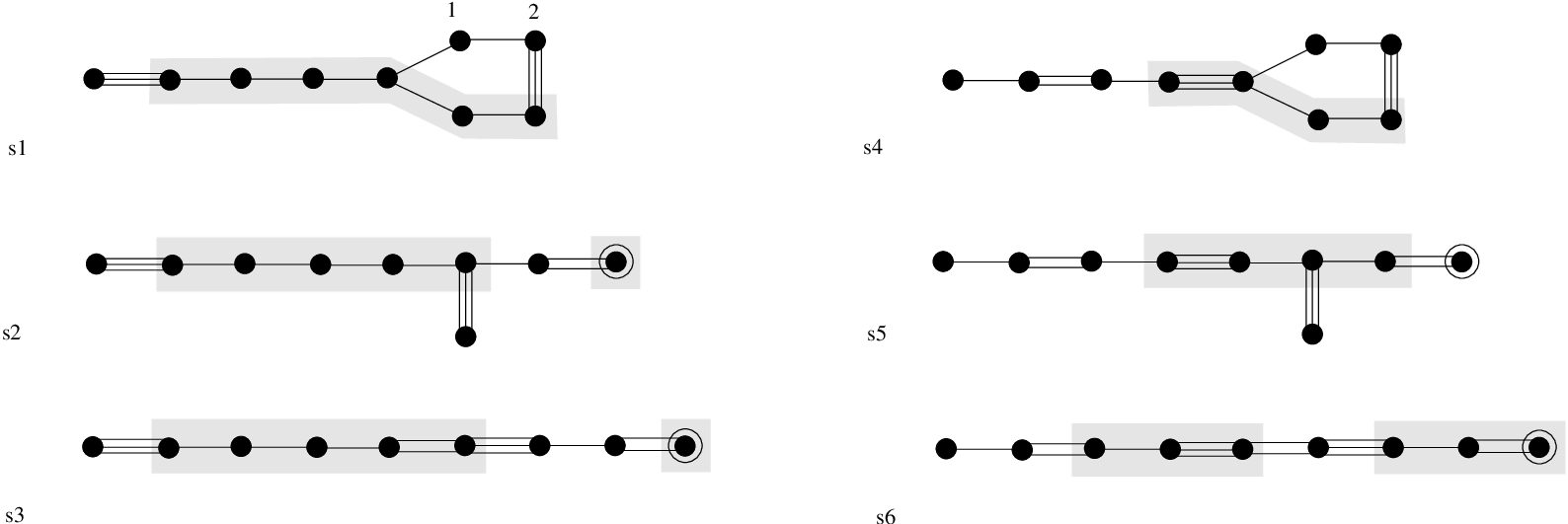,width=0.75\linewidth}
\end{center}
\end{table}

First, consider the diagram $\Sigma_1$. Let $t_1$ and $t_2$ be the nodes of $\Sigma_1$ marked in Table~\ref{7_f4_h4}.
Without loss of generality we may assume that neither $t_1$ nor $t_2$ is an end of the dotted edge
(here we use the symmetry of the diagram $\Sigma_1$). Let $S_2\subset \Sigma_1$ be a diagram of the type $A_6$
that does not contain  the nodes $t_1$ and $t_2$. Then $\Sigma(P)$ contains a diagram from the list
$L'(\l S_2,t_1,t_2\r,5,7)$, which is empty.

For the diagrams   $\Sigma_2,\dots,\Sigma_6$ denote by $S_2$ a subdiagram marked by a gray box.
In cases of $\Sigma_4$ and $\Sigma_5$ the diagram $S_2$ is of order 4, and it has only 2 good neighbors 
(or non-neighbors) in $\Sigma_i$,
so it has at least 2 more good neighbors (or non-neighbors) in $\Sigma(P)$, one of which is joined with  $\Sigma_i$ 
without dotted edges.
In cases of $\Sigma_2$,  $\Sigma_3$, and $\Sigma_6$, the diagram $S_2$ is of order 6, and it has only 1  
good neighbor (or non-neighbor) in $\Sigma_i$,
so, it has another one in  $\Sigma(P)\setminus \Sigma_i$ (and this good neighbor or non-neighbor 
cannot be joined with $\Sigma_i$
by a dotted edge since $S_2$ contains an end of the dotted edge).
Therefore, $\Sigma(P)$ contains a diagram from the list $L(\Sigma_i,5,7,S_2)$, where $i=2,\dots,6$.
For $i=2,3,4$ the lists are empty. For $i=5$ and $i=6$ the lists consist of a unique entry $\Sigma^{56}$ shown in
Fig.~\ref{7_f4h4_2}. Denote by $S_3\subset \Sigma^{56}$ a subdiagram of order 6 marked by a gray box.
It has only one good neighbor (and no non-neighbors) in $\Sigma^{56}$ and contains an end of the dotted edge.
Hence,  $\Sigma(P)$ contains a diagram from the list $L(\Sigma^{56},5,7,S_3)$, which is empty.

\begin{figure}[!h]
\begin{center}
\psfrag{s}{$\Sigma^{56}$}
\epsfig{file=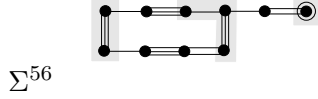,width=0.327\linewidth}
\caption{Treating the diagrams $\Sigma_5$ and $\Sigma_6$, see Lemma~\ref{7_h4,f4}.}
\label{7_f4h4_2}
\end{center}
\end{figure}

\end{proof}

\begin{lemma}
\label{7_h3}
$\Sigma(P)$  contains no subdiagram of the type $H_3$.

\end{lemma}

\begin{proof}
Suppose that $S_0\subset \Sigma(P)$ is a subdiagram of the type $H_3$.
Then $P(S_0)$ is a 4-polytope whose Coxeter diagram contains at most 1 dotted edge,
so it is either a simplex, or an Esselmann polytope, or a 4-prism, or a 4-polytope with 7 facets.
Since $\o S_0=\Sigma_{S_0}$ contains neither multi-multiple edges nor subdiagrams of the types $H_4$ and $F_4$,
we are left with only three possibilities for $\o S_0$ shown in Fig.~\ref{6_h2_}(a).
For each of these diagrams consider a subdiagram $\Sigma'$ of order 5 shown in Fig.~\ref{7_h3_},
and let $S_1\subset \Sigma'$ be a subdiagram of order 4 marked by a gray block.
Notice that $S_1$ has at least one good neighbor or non-neighbor in $\Sigma(P)\setminus \l S_0,\o S_0\r$,
so $\Sigma(P)$ contains a diagram from the list $L'(\Sigma',5,7,S_1)$, where $\Sigma'$ ranges over the three diagrams 
shown in Fig.~\ref{7_h3_}.
These lists are empty, and the lemma is proved.

\begin{figure}[!h]
\begin{center}
\psfrag{2,3}{\scriptsize $2,3$}
\epsfig{file=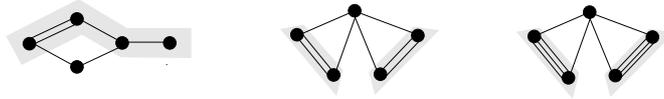,width=0.7\linewidth}
\caption{To the proof of Lemma~\ref{7_h3}.}
\label{7_h3_}
\end{center}
\end{figure}

\end{proof}

\begin{lemma}
\label{7_h2}
$\Sigma(P)$  contains no subdiagram of the type $G_2^{(5)}$.

\end{lemma}

\begin{proof}
Suppose that $S_0\subset \Sigma(P)$ is a subdiagram of the type $G_2^{(5)}$.
Then $P(S_0)$ is a 5-polytope with at most one pair of non-intersecting facets.
By the Main Theorem in dimension 5, this implies that $P(S_0)$ has at most 8 facets.
However, any diagram of a 5-polytope with at most 8 facets contains either 2 dotted edges or 
a subdiagram of the types $H_4$ or $F_4$.
Together with Lemma~\ref{7_h4,f4} this proves the lemma.   

\end{proof}

Applying Lemmas~\ref{b_d}, and~\ref{b_2}, we obtain the following result.

\begin{lemma}
\label{d7}
The Main Theorem holds in dimension 7. 

\end{lemma}

\subsection{Dimension 8}

Let $P$ be an 8-dimensional compact hyperbolic Coxeter polytope such that $\Sigma(P)$ contains a unique dotted edge and 
$P$ has at least 12 facets.

\begin{lemma}
\label{8_mult}
$\Sigma(P)$  contains no multi-multiple edges.

\end{lemma}

\begin{proof}
Suppose that $S_0\subset \Sigma(P)$ is a multi-multiple edge of the maximum multiplicity in $\Sigma(P)$. 
Then $S_0$ has no good neighbors and
$P(S_0)$ is a Coxeter 6-polytope with at most 1 pair of non-intersecting facets.
Since the Main Theorem is already proved in dimension 6, this implies that
$P(S_0)$ has at most $9$ facets and
$\o S_0$ is one of the 3 diagrams $\Sigma_1$, $\Sigma_2$, $\Sigma_3$ shown in Fig.~\ref{8_multi}.

\begin{figure}[!h]
\begin{center}
\psfrag{a}{(a)}
\psfrag{b}{(b)}
\psfrag{s1}{\small $\Sigma_1$}
\psfrag{s2}{\small $\Sigma_2$}
\psfrag{s3}{\small $\Sigma_3$}
\psfrag{8}{\small $8$}
\psfrag{10}{\scriptsize $10$}
\epsfig{file=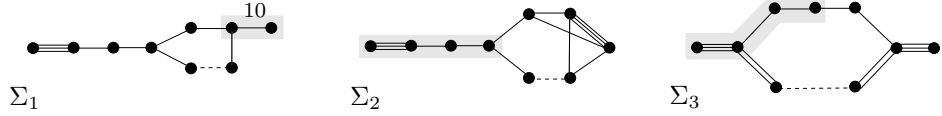,width=0.97\linewidth}
\caption{To the proof of Lemma~\ref{8_mult}.}
\label{8_multi}
\end{center}
\end{figure}

Consider the diagram $\Sigma_1$. It contains a subdiagram $S_1$ of the type $G_2^{(10)}$ such that 
$\o S_1=\Sigma_{S_1}$ contains no dotted edge. Since $P(S_1)$ is a 6-polytope, this is impossible.

Consider the diagram $\Sigma_2$. It contains a subdiagram $S_1$ of the type $H_4$ (marked by a gray box) such that 
$\o S_1=\Sigma_{S_1}$ contains no dotted edge. $P(S_1)$ is a 4-polytope,
so $\o S_1$ is either a Lann\'er diagram of order 5 or an Esselmann diagram. At the same time,
$\o S_1$ contains a multi-multiple edge $S_0$ and a Lann\'er diagram of order 3 with one triple edge and two simple edges.
This is impossible for an Esselmann diagram as well as for a Lann\'er diagram of order 5.

Consider the diagram $\Sigma_3$. It contains a subdiagram $S_1$ of the type $H_4$ such that 
$\o S_1=\Sigma_{S_1}$ contains no dotted edge. At the same time, $\o S_1$ contains 
a multi-multiple edge $S_0$ and a Lann\'er diagram of order 3 with one triple edge, one double edge, and one empty edge.
This is possible only if $\o S_1$ is an Esselmann diagram and $S_0=G_2^{(10)}$. 
In particular, this implies that any multi-multiple edge in 
$\Sigma(P)$ is of the type $G_2^{(10)}$. 
Denote by $\Sigma_3'$ the diagram $\Sigma_3$ with one end of the dotted edge discarded.
Let $S_2\subset \Sigma_3'$ be a subdiagram of the type $B_6$.
It has only two non-neighbors (and no good neighbors) in $\l  S_0, \Sigma_3\r$, so there exists either a good neighbor or a non-neighbor 
$x$ of $S_2$, such that $x\notin \l  S_0, \Sigma_3\r$ and the diagram $\l x, S_0, \Sigma_3'\r$ contains no dotted edges.
Since any multi-multiple edge in $\Sigma(P)$ is of the type $G_2^{(10)}$, the number of such diagrams is finite. 
None of these diagrams has zero determinant, so the lemma is proved.

\end{proof}

\begin{lemma}
\label{8_h4,f4}
$\Sigma(P)$  contains no subdiagrams of the types $H_4$ and $F_4$.

\end{lemma}

\begin{figure}[!h]
\begin{center}
\psfrag{a}{(a)}
\psfrag{b}{(b)}
\psfrag{6}{\small $6$}
\psfrag{8}{\small $8$}
\psfrag{3,4}{\tiny $3,4$}
\epsfig{file=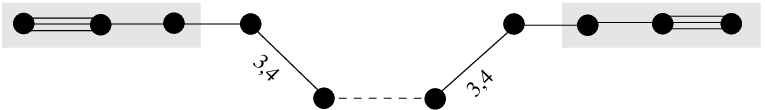,width=0.37\linewidth}
\caption{The diagram  $\l S_0,y_2,\o  S_0 \r$, see Lemma~\ref{8_h4,f4}. }
\label{8_h4_2}
\end{center}
\end{figure}

\begin{proof}
Suppose that $S_0\subset \Sigma(P)$ is a subdiagram of the type $H_4$ or $F_4$.
$S_0$ has no good neighbors, so $\Sigma(P)$ contains a diagram from the list 
$L_1(H_4,8)$ or $L_1(F_4,8)$. 
The union of these lists consists of 9 diagrams $\Sigma_1,\dots,\Sigma_9$, see Table~\ref{8_h4_}.
One can note that for any of diagrams  $\Sigma_1,\dots,\Sigma_9$ the diagram $\o S_0$ is a linear Lann\'er subdiagram 
containing a subdiagram of the type $H_4$, and $\o S_0\subset \Sigma_i$
(by a linear diagram we mean a connected diagram without nodes of valency greater than 2).
This implies that we can always start from the diagram $S_0$ of the type $H_4$, so $\Sigma(P)$ must contain 
one of the diagrams $\Sigma_1,\dots,\Sigma_6$, and we do not need to consider   
the diagrams $\Sigma_7$, $\Sigma_8$, and  $\Sigma_9$. Moreover, notice that $y_1$ (which is a unique bad neighbor of 
$S_0$ in $\Sigma_i$) is always a bad neighbor of a unique subdiagram $S_2 \subset \o S_0$ of the type $H_4$.
By construction (see Lemma~\ref{l1}), this implies that there exists a non-neighbor $y_2\notin\Sigma_i$ 
of $S_2$ joined with $\o S_0 \setminus S_2$ by a dotted edge. 
Starting from $S_2$ instead of $S_0$, we obtain (by symmetry) that $\o S_2$ is also a linear Lann\'er diagram of order 5. 
Since $\l S_0,y_2\r\subset \o S_2$, we see that both $\l S_0,y_2 \r$ and  $ \o  S_0 $ are  linear Lann\'er diagrams, 
and $y_2$ is joined with $\o S_0\setminus S_2$ by a dotted edge.
Thus, we obtain three possibilities for the subdiagram $\l S_0,y_2,\o  S_0 \r$, see Fig.~\ref{8_h4_2}.
For each of these 3 diagrams we solve the equation $\det (\l S_0,y_2,\o  S_0 \r)=0$ and find the weight of the dotted edge.
Consider a diagram $S_3\subset  \l S_0,y_2,\o  S_0 \r$ of the type $H_3+H_3$ (it is marked on Fig.~\ref{8_h4_2}). 
$S_3$ has four good neighbors and non-neighbors in total in $ \l S_0,y_2,\o  S_0 \r$, 
while $\o S_3$ has at least three dotted edges (one coming from a dotted edge of $\Sigma(P)$ 
and two coming from simple or double edges). This implies that $S_3$ has at least one good neighbor or a non-neighbor in 
$\Sigma(P)\setminus  \l S_0,y_2,\o  S_0 \r$. So, $\Sigma(P)$ contains a diagram from the list $L'(\l S_0,y_2,\o  S_0 \r,5,8)$.
This list consists of a unique diagram, which is a diagram of a Coxeter 8-polytope with 11 facets (see Fig.~\ref{10_}). 
By~\cite[Lemma~1]{nodots}, this diagram cannot be a subdiagram of $\Sigma(P)$. 

\end{proof}

\begin{table}[!h]
\begin{center}
\caption{The lists $L_1(H_4,8)$ and $L_1(F_4,8)$. }
\label{8_h4_}
\vspace{10pt}
\psfrag{a}{(a)}
\psfrag{b}{(b)}
\psfrag{6}{\small $6$}
\psfrag{8}{\small $8$}
\psfrag{10}{\small $10$}
\psfrag{s1}{\small $\Sigma_1$}
\psfrag{s2}{\small $\Sigma_2$}
\psfrag{s3}{\small $\Sigma_3$}
\psfrag{s4}{\small $\Sigma_4$}
\psfrag{s5}{\small $\Sigma_5$}
\psfrag{s6}{\small $\Sigma_6$}
\psfrag{s7}{\small $\Sigma_7$}
\psfrag{s8}{\small $\Sigma_8$}
\psfrag{s9}{\small $\Sigma_9$}
\epsfig{file=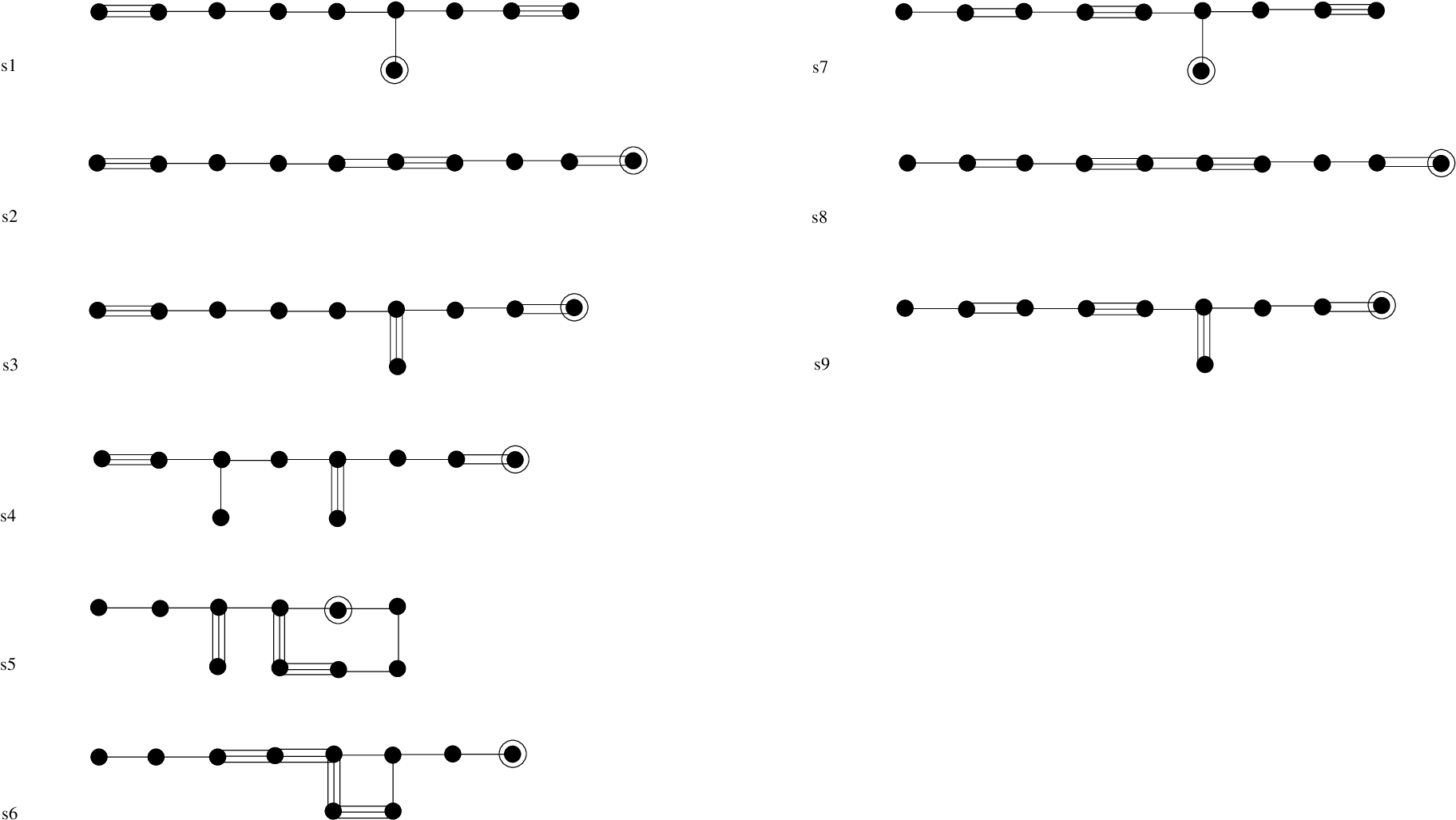,width=0.9\linewidth}
\end{center}
\end{table}

\begin{lemma}
\label{8_h3}
$\Sigma(P)$  contains no subdiagram of the type $H_3$.

\end{lemma}

\begin{proof}
Suppose that $S_0\subset \Sigma(P)$ is a subdiagram of the type $H_3$.
Then $P(S_0)$ is a 5-polytope with at most one pair of non-intersecting facets.
By the Main Theorem in dimension 5, this implies that $P(S_0)$ has at most 8 facets.
However, any diagram of a 5-polytope with at most 8 facets either contains 2 dotted edges or 
contains a subdiagram of the types $H_4$ or $F_4$.
Together with Lemma~\ref{8_h4,f4}, this proves the lemma.

\end{proof}

\begin{lemma}
\label{8_h2}
$\Sigma(P)$  contains no subdiagram of the type $G_2^{(5)}$.

\end{lemma}

\begin{proof}
Suppose that $S_0\subset \Sigma(P)$ is a subdiagram of the type $G_2^{(5)}$.
Then $P(S_0)$ is a 6-polytope with at most one pair of non-intersecting facets.
By the Main Theorem in dimension 6, this implies that $P(S_0)$ has at most 9 facets,
so $P(S_0)$ has exactly 9 facets.
However, any diagrams of a 6-polytope with 9 facets 
contains a subdiagram of the type $H_4$.
Together with Lemma~\ref{8_h4,f4}, this proves the lemma.   

\end{proof}

\noindent
As in dimensions 6 and 7, we apply Lemmas~\ref{b_d} and~\ref{b_2} to obtain  

\begin{lemma}
\label{d8}
The Main Theorem holds in dimension 8. 

\end{lemma}

\subsection{Dimension 9}

\begin{lemma}
\label{9}
The Main Theorem holds in dimension 9.

\end{lemma}

\begin{proof}
Suppose that the lemma is broken.
Let $P$ be a 9-dimensional compact hyperbolic Coxeter polytope such that $\Sigma(P)$ contains a unique dotted edge and 
$P$ has at least 13 facets.

\smallskip
\noindent
$\bullet$ 
{\it $\Sigma(P)$ contains no multi-multiple edges.}\\
Indeed, if $S_0\subset \Sigma(P)$ is a multi-multiple edge, then $P(S_0)$ is a 7-polytope with at most one pair of 
non-intersecting facets, so $P(S_0)$ is a 7-polytope with at most 10 facets, which does not exists.

\smallskip
\noindent
$\bullet$ 
{\it $\Sigma(P)$ contains no subdiagrams of the types $H_4$ and $F_4$.}\\
Suppose that $S_0\subset \Sigma(P)$ is a subdiagram of the type $H_4$ or $F_4$.
Then $P(S_0)$ is a 5-polytope with at most one pair of non-intersecting facets,
so $P(S_0)$ is a 5-polytope with at most 8 facets.
Since $\o S_0=\Sigma_{S_0}$ contains no multi-multiple edges and at most one dotted edge,
there are only three possibilities for the diagram $\o S_0$, see Fig.~\ref{9_}(a)--(c).
For each of these cases we choose a subdiagram $\Sigma_1$ of order 6 shown in Fig.~\ref{9_}(d)--(f) respectively,
and denote by $S_1\subset \Sigma_1$ a subdiagram of the type $H_4$ or $F_4$ marked by a gray box.
Let $S_2\subset S_0$ be a subdiagram of the type $H_3$ or $B_3$ (if $S_0$ is of the type $H_4$ or 
$F_4$, respectively). Let $S_3=\l S_1,S_2 \r$.
Notice that $S_3$ has 3 good neighbors and non-neighbors in total in $\l S_0,\o S_0 \r$,
two of which are the ends of the dotted edge. Hence, by Lemma~\ref{polygon}, $S_3$ has 
at least one good neighbor or non-neighbor in $\Sigma(P)\setminus \l S_0,\o S_0 \r$.
Therefore, $\Sigma(P)$ contains a diagram from the list $L'(\l S_0,\Sigma_1 \r,5,9,S_3)$,
where $\Sigma_1$ ranges over the diagrams shown in Fig.~\ref{9_}(d)--(f).
The lists are empty, and the statement is proved.

\begin{figure}[!h]
\begin{center}
\psfrag{a}{\small (a)}
\psfrag{b}{\small (b)}
\psfrag{c}{\small (c)}
\psfrag{d}{\small (d)}
\psfrag{e}{\small (e)}
\psfrag{f}{\small (f)}
\psfrag{6}{\small $6$}
\psfrag{8}{\small $8$}
\psfrag{10}{\small $10$}
\psfrag{3,4}{\tiny $3,\!4$}
\psfrag{2,3,4}{\tiny $2,\!3,\!4$}
\psfrag{2,3}{\tiny $2,\!3$}
\epsfig{file=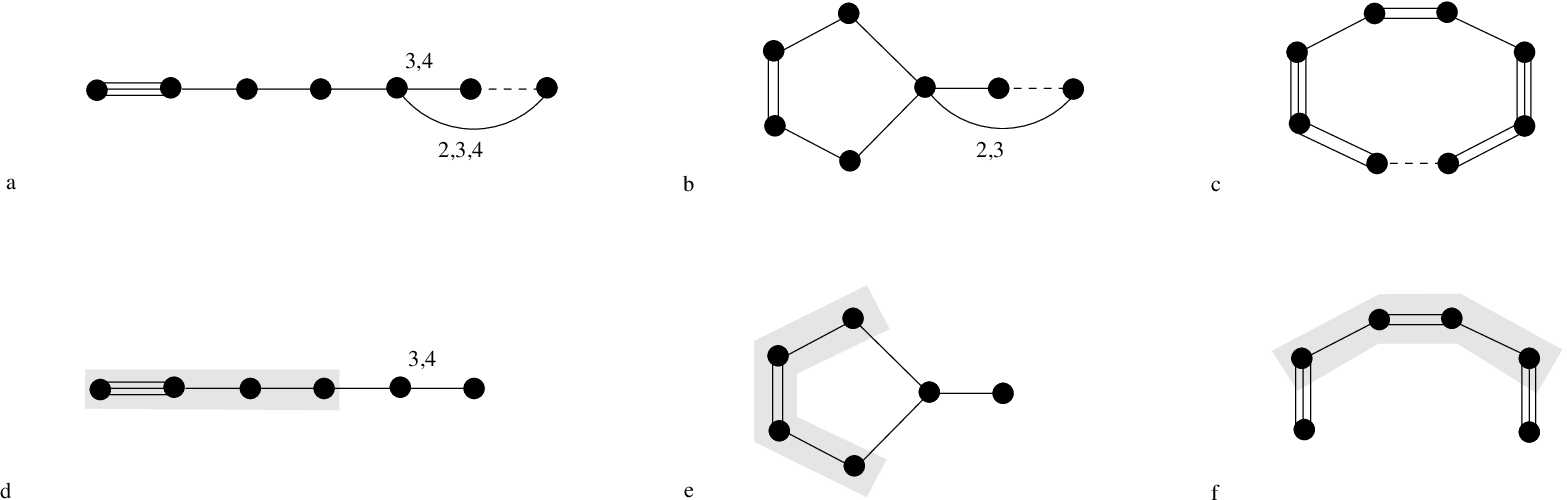,width=0.8\linewidth}
\caption{To the proof of Lemma~\ref{9}.}
\label{9_}
\end{center}
\end{figure}

\smallskip
\noindent
$\bullet$ 
{\it $\Sigma(P)$ contains no subdiagrams of the types $H_3$.}\\
Indeed, if $S_0\subset \Sigma(P)$ is a subdiagram of the type $H_3$, 
then $P(S_0)$ is a 6-polytope with at most one pair of non-intersecting facets.
However, a diagram of any such a polytope contains a subdiagram of the type $H_4$.

\smallskip
\noindent
$\bullet$ 
{\it $\Sigma(P)$ contains no subdiagrams of the types $G_2^{(5)}$.}\\
If $S_0\subset \Sigma(P)$ is a subdiagram of the type $G_2^{(5)}$,
 then $P(S_0)$ is a 7-polytope with at most one pair of non-intersecting facets, which does not exists.

Now, we apply Lemmas~\ref{b_d} and~\ref{b_2}, which finishes the proof. 

\end{proof}

\subsection{Dimension 10}

\begin{lemma}
\label{10}
The Main Theorem holds in dimension 10.

\end{lemma}

\begin{proof}
Suppose that the lemma is broken.
Let $P$ be a 10-dimensional compact hyperbolic Coxeter polytope such that $\Sigma(P)$ contains a unique dotted edge. 

\smallskip
\noindent
$\bullet$ 
{\it $\Sigma(P)$ contains no multi-multiple edges.}\\
Indeed, if $S_0\subset \Sigma(P)$ is a multi-multiple edge, then $P(S_0)$ is a 8-polytope with at most one pair of 
non-intersecting facets,
so  $P(S_0)$ is a 8-polytope with at most 11 facets.
There exists a unique such a polytope, its diagram is shown in Fig.~\ref{10_}.
Let $S_1\subset \o S_0$ be a subdiagram of the type $H_4$. Then $\o S_1$ contains no dotted edges, and $P(S_1)$
is a Coxeter 6-polytope with mutually intersecting facets, which is impossible.  

\begin{figure}[!h]
\begin{center}
\psfrag{a}{(a)}
\psfrag{b}{(b)}
\psfrag{6}{\small $6$}
\psfrag{8}{\small $8$}
\psfrag{10}{\small $10$}
\epsfig{file=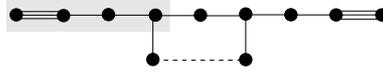,width=0.4\linewidth}
\caption{A unique 8-polytope with 11 facets.}
\label{10_}
\end{center}
\end{figure}

\smallskip
\noindent
$\bullet$ 
{\it $\Sigma(P)$ contains no subdiagrams of the types $H_4$ and $F_4$.}\\
Suppose that $S_0\subset \Sigma(P)$ is a subdiagram of the type $H_4$ or $F_4$.
Then $P(S_0)$ is a 6-polytope with at most one pair of non-intersecting facets,
so $P(S_0)$ is a 6-polytope with exactly  9 facets.
There are 3 such polytopes (see Fig.~\ref{8_multi}), each contains a subdiagram $S_1$ of the type $H_4$
such that $\o S_1$ contains no dotted edges.
So, $P(S_1)$ is a 6-polytope with mutually intersecting facets, which is impossible.  

\smallskip
\noindent
$\bullet$ 
{\it $\Sigma(P)$ contains no subdiagrams of the types $H_3$.}\\
Indeed, if $S_0\subset \Sigma(P)$ is a subdiagram of the type $H_3$, 
then $P(S_0)$ is a 7-polytope with at most one pair of non-intersecting facets.
This implies that $P(S_0)$ is a 7-polytope with at most 10 facets, which is impossible.

\smallskip
\noindent
$\bullet$ 
{\it $\Sigma(P)$ contains no subdiagrams of the types $G_2^{(5)}$.}\\
As it was already shown, the diagram of the type  $G_2^{(5)}$ cannot have good neighbors,
so the proof coincides with the reasoning used for multi-multiple edges.

Applying Lemmas~\ref{b_d} and~\ref{b_2}, we complete the proof. 

\end{proof}

\subsection{Dimension 11}

\begin{lemma}
\label{11}
The Main Theorem holds in dimension 11.

\end{lemma}

\begin{proof}
Suppose that the lemma is broken.
Let $P$ be a 11-dimensional compact hyperbolic Coxeter polytope such that $\Sigma(P)$ contains a unique dotted edge. 

\smallskip
\noindent
$\bullet$ 
{\it $\Sigma(P)$ contains no multi-multiple edges.}\\
If $S_0\subset \Sigma(P)$ is a multi-multiple edge, then $P(S_0)$ is a 9-polytope with at most one pair of 
non-intersecting facets.

\smallskip
\noindent
$\bullet$ 
{\it $\Sigma(P)$ contains no subdiagrams of the types $H_4$ and $F_4$.}\\
Indeed, if $S_0\subset \Sigma(P)$ is a  subdiagram of the type $H_4$ or $F_4$, then $P(S_0)$ is a 7-polytope 
with at most one pair of non-intersecting facets,
which is impossible.

\smallskip
\noindent
$\bullet$ 
{\it $\Sigma(P)$ contains no subdiagrams of the types $H_3$.}\\
If $S_0\subset \Sigma(P)$ is a subdiagram of the type $H_3$, 
then $P(S_0)$ is a 8-polytope with at most one pair of non-intersecting facets.
However, the diagram of a unique such a polytope contains a subdiagram of the type $H_4$.

\smallskip
\noindent
$\bullet$ 
{\it $\Sigma(P)$ contains no subdiagrams of the types $G_2^{(5)}$.}\\
Again, we follow the proof for multi-multiple edges. 

Application of Lemmas~\ref{b_d} and~\ref{b_2} finishes the proof. 

\end{proof}

\subsection{Dimension 12}

\begin{lemma}
\label{12}
The Main Theorem holds in dimension 12.

\end{lemma}

\begin{proof}
Suppose that the lemma is broken.
Let $P$ be a 12-dimensional hyperbolic Coxeter polytope such that $\Sigma(P)$ contains a unique dotted edge. 

\smallskip
\noindent
$\bullet$ 
{\it $\Sigma(P)$ contains no subdiagrams of the types $H_4$ and $F_4$.}\\
Indeed, if $S_0\subset \Sigma(P)$ is a  subdiagram of the type $H_4$ or $F_4$, then $P(S_0)$ is a 8-polytope with at most 
one pair of non-intersecting facets.
So, $\o S_0$ is the diagram shown in Fig.~\ref{10_}. 
However, the latter diagram contains a subdiagram $S_1$
of the type $H_4$  such that $\o S_1$ contains no dotted edges, which is impossible.

\smallskip
\noindent
$\bullet$ 
{\it $\Sigma(P)$ contains no subdiagrams of the types $H_3$ and  $G_2^{(k)}$, $k\ge 5$.}\\
If $S_0\subset \Sigma(P)$ is a subdiagram of the type $H_3$ or $G_2^{(k)}$, $k\ge 5$,
then $P(S_0)$ is a $d$-polytope with at most one pair of non-intersecting facets, where $d=9$ or 10, 
which is impossible.

Again, we complete the proof applying Lemmas~\ref{b_d} and~\ref{b_2}. 

\end{proof}

\subsection{Large dimensions}
To complete the proof of the Main Theorem, we prove the following lemma.

\begin{lemma}
The Main Theorem holds in dimensions $d>12$.

\end{lemma}

\begin{proof}
Suppose that the lemma is broken, and let $P$ be a $d$-dimensional compact 
hyperbolic Coxeter polytope such that $\Sigma(P)$ contains a unique dotted edge
($d>12$). We may assume that the Main Theorem holds in all dimensions less than $d$. 
Suppose that   $\Sigma(P)$ contains a subdiagram $S_0$ of the type $H_4$ or $F_4$.
Then  $P(S_0)$ is a $d$-polytope with at most one pair of non-intersecting facets, where $d\ge 9$, 
which is impossible.
Similarly, $\Sigma(P)$ contains no subdiagrams of the types $H_3$ and $G_2^{(k)}$, $k\ge 5$.

As usual, Lemmas~\ref{b_d} and~\ref{b_2} imply that such a polytope $P$ does not exist. 

\end{proof}

\section*{Appendix}
In this section we list all compact hyperbolic Coxeter polytopes with exactly one pair of non-intersecting facets.
Table~\ref{answer1} contains Coxeter diagrams of $d$-polytopes with $d+3$ facets, the list is reproduced from~\cite{d+3}. 

\begin{table}
\begin{center}
\caption{Compact hyperbolic Coxeter $d$-polytopes with $d+3$ facets and exactly one pair of non-intersecting facets} 
\label{answer1}
\end{center}
\underline{\bf d=4}

\begin{tabular}{ccp{0.05\linewidth}cp{0.05\linewidth}c}
&
\psfrag{8}{\small $8$}
\epsfig{file=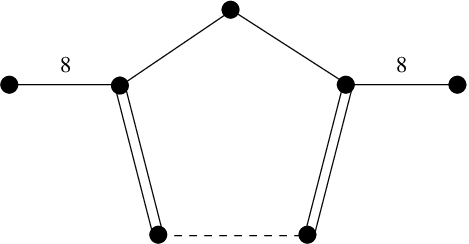,width=0.32\linewidth}&&
\psfrag{8}{\small $8$}
\epsfig{file=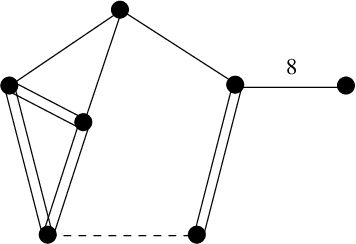,width=0.24\linewidth}&&
\epsfig{file=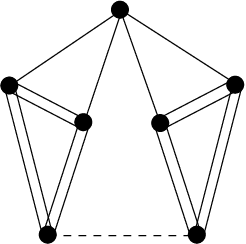,width=0.165\linewidth}\\
\\
&
\psfrag{10}{\small $10$}
\epsfig{file=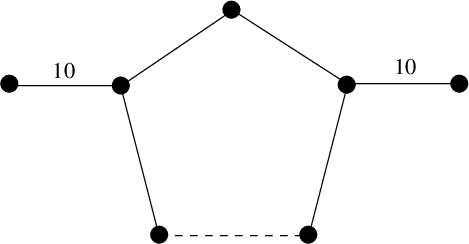,width=0.32\linewidth}&&
\psfrag{10}{\small $10$}
\epsfig{file=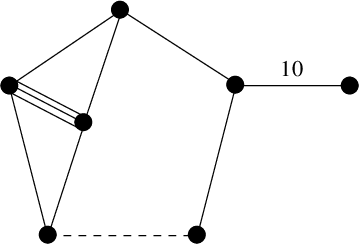,width=0.24\linewidth}&&
\epsfig{file=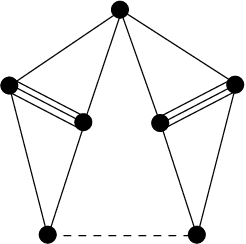,width=0.165\linewidth}\\
\\
&\multicolumn{5}{c}{
\begin{tabular}{p{0.05\linewidth}cp{0.1\linewidth}c}
&
\psfrag{8}{\small $8$}
\epsfig{file=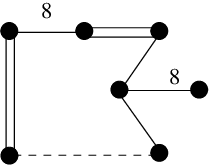,width=0.17\linewidth}&&
\psfrag{8}{\small $8$}
\epsfig{file=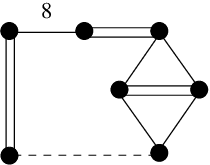,width=0.17\linewidth}\\
\\
\end{tabular}
}\\

\end{tabular}

\noindent
\underline{\bf d=5}

\begin{tabular}{p{0.08\linewidth}c}
&\epsfig{file=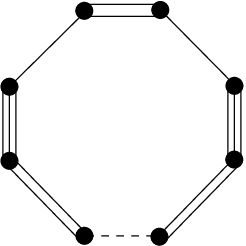,width=0.175\linewidth}\\
\\
\end{tabular}

\noindent
\underline{\bf d=6}

\begin{tabular}{p{0.01\linewidth}cp{0.05\linewidth}cp{0.05\linewidth}c}
&
\raisebox{0.08\linewidth}{\epsfig{file=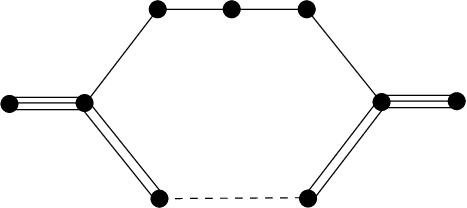,width=0.32\linewidth}}&&
\psfrag{10}{\small $10$}
\raisebox{0.08\linewidth}{\epsfig{file=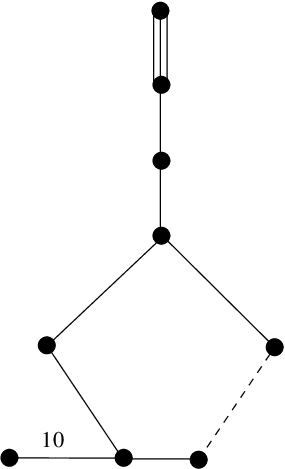,width=0.18\linewidth}}&&
\psfrag{10}{\small $10$}
\raisebox{0.04\linewidth}{\epsfig{file=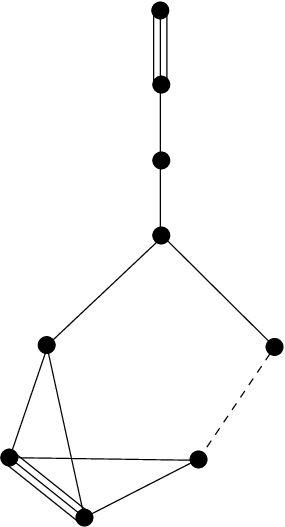,width=0.18\linewidth}}
\end{tabular}

\vspace{-0.03\linewidth}

\noindent
\underline{\bf d=8}

\begin{tabular}{p{0.08\linewidth}c}
&\epsfig{file=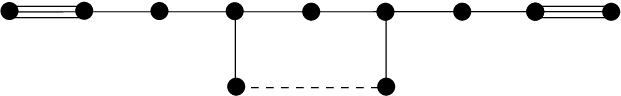,width=0.50\linewidth}
\end{tabular}
\end{table}

\end{document}